\documentclass[11pt]{article}

\usepackage{epsfig,lscape}
\usepackage{amsmath,amssymb}
\usepackage{longtable}
\usepackage{algorithm,algorithmic}
\usepackage{color,varioref}
\usepackage{verbatim,graphicx}
\usepackage{amsthm,amsfonts,latexsym,rawfonts}
\usepackage{hyperref}
\usepackage{url}
\RequirePackage{ifpdf}
\usepackage{comment}
\usepackage{amsthm}
\usepackage{amsmath}
\usepackage{amssymb}
\usepackage{cases}
\usepackage{latexsym}
\usepackage{graphicx}
\usepackage{fancybox,color}
\usepackage{longtable}
\usepackage{color}
\usepackage{lscape}
\usepackage{graphicx}
\usepackage{footnote}
\makesavenoteenv{tabular}
\usepackage{float}
\definecolor{lred}{rgb}{1,0.8,0.8}
\definecolor{lblue}{rgb}{0.8,0.8,1}
\definecolor{dred}{rgb}{0.6,0,0}
\definecolor{dblue}{rgb}{0,0,0.5}
\definecolor{mygreen}{RGB}{28,172,0} 
\definecolor{grey}{rgb}{0.98,0.98,0.98}
\definecolor{mylilas}{RGB}{170,55,241}
\usepackage{listings}
\usepackage{url}

\def\cA{{\cal A}} 
\def\cB{{\cal B}} \def\cE{{\cal E}} 

\def\cD{{\cal D}}
\def\cU{{\cal U}}
\def\cM{{\cal M}}
\def\cP{{\cal P}}

\def\cH{{\cal H}}

 \def\S{{\cal S}}

\newcommand{\inprod}[2]{\langle #1 , #2 \rangle}

\def\norm#1{\|#1\|}
\def\abs#1{\left| #1\right|}

\def\cM{{\cal M}}
\def\cQ{{\cal Q}}
\def\cK{{\cal K}}
\def\cI{{\cal I}}

\def\inprod#1#2{\langle#1, \, #2\rangle}

\def\diag{{\rm diag}}

\def\sig{\sigma}

\def\cX{{\cal X}}

\newcommand{\T}{\mbox{\textrm{\tiny{T}}}}

\def\mc{\multicolumn}

\def\by{\bar{y}}

\def\inprod#1#2{\langle #1,\,#2\rangle}
\def\norm#1{\|#1\|}

\def\mc{\multicolumn}

\def\mc{\multicolumn}

\def\by{\bar{y}}

\def\inprod#1#2{\langle #1,\,#2\rangle}
\def\norm#1{\|#1\|}

\def\SDPNALP{{\mbox{{\sc Sdpnal}$+$}} }
\def\SDPNALPP{{\mbox{{\sc Sdpnal}$+$. }}}
\renewcommand{\Re}{\mathbb{R}}
\def\S{\mathbb{S}}

\begin{document}

\title{\bf SDPNAL$+$: A {\sc Matlab} software for semidefinite programming with bound constraints (version 1.0)
	}
\author{
	Defeng Sun\thanks{Department of Applied Mathematics, The Hong Kong Polytechnic University, Hung Hom,
	Hong Kong ({\tt defeng.sun@polyu.edu.hk}).
The research of this author is partially supported by a start-up research grant from the Hong Kong Polytechnic University.},
\;
Kim-Chuan Toh\thanks{Department of Mathematics, and Institute of Operations Research and Analytics, National University of Singapore, 10 Lower Kent Ridge Road, Singapore
		({\tt mattohkc@nus.edu.sg}).
		The research of this author is supported in part by the Ministry of
Education, Singapore, Academic Research Fund under Grant R-146-000-256-114.
		},
\;
Yancheng Yuan\thanks{Department of Mathematics, National University of Singapore, 10 Lower Kent Ridge Road, Singapore ({\tt yuanyancheng@u.nus.edu}).},
\;
Xin-Yuan Zhao\thanks{Beijing Institute for Scientific and Engineering Computing,
Beijing University of Technology, 100 Pingleyuan,
Chaoyang District, Beijing 100124, People's Republic of China ({\tt xyzhao@bjut.edu.cn}). The research of this author was  supported by the National Natural Science Foundation of China under projects No.11871002 and the General Program of Science and Technology of Beijing Municipal Education Commission.}
}
\date{December 28, 2018}
\maketitle

\textbf{Abstract.}
\SDPNALP is a {\sc Matlab} software package that implements an
augmented Lagrangian based method to solve large scale semidefinite programming problems with bound constraints.  The implementation was initially based
on a majorized semismooth Newton-CG augmented Lagrangian method,  here
we designed it within
an inexact symmetric Gauss-Seidel based semi-proximal
ADMM/ALM (alternating direction method of multipliers/augmented Lagrangian method)
framework
for  the purpose of deriving simpler stopping conditions and
closing the gap between the practical implementation of the algorithm and the
theoretical algorithm.
The basic code is written in {\sc Matlab}, but some subroutines in C language are incorporated via Mex files.
We also design a convenient interface for  users to input their SDP models
into the solver.
Numerous problems arising from combinatorial optimization and binary integer quadratic programming problems have been tested to evaluate the performance of the solver.
Extensive numerical experiments conducted in
[Yang, Sun, and  Toh, Mathematical Programming Computation, 7 (2015), pp. 331--366]
 show that the proposed method is quite efficient and robust,
 in that it is able to solve 98.9\% of the 745  test instances of SDP problems
 arising from various applications to the accuracy of
 $ 10^{-6}$ in the relative KKT residual.

\medskip
\noindent \textbf{Keywords: } Semidefinite programming, Augmented Lagrangian, Semismooth Newton-CG method, Matlab software package.

\section{Introduction}


Let $\S^n$ be the space of $n \times n$ real symmetric matrices and $\S^n_+$ be the cone of  positive semidefinite matrices in $\S^n$.
For any $X \in \S^n$, we may sometimes write $X \succeq 0$  to indicate that
$X\in \S^n_+$.
Let $\cP = \{X \in \S^n: L \leq X \leq U\}$, where $L, U$ are given $n\times n$ symmetric matrices whose
elements are allowed to take the values $-\infty$ and $+\infty$, respectively.
Consider the semidefinite programming (SDP) {problem}:
\begin{eqnarray*} \label{eq-sdp-org}
({\rm \bf SDP}) \quad  \min \Big\{ \inprod{C}{X} \mid \cA(X) = b,\;
 l \leq \cB(X) \leq u, \;
X \in \S^n_+, \; X\in \cP\Big\},
\end{eqnarray*}
where $b\in\Re^m$,  and $C \in \S^n$ are given data, $\cA:\S^n\rightarrow \Re^m$ and $\cB:\S^n\rightarrow \Re^{p}$ are two given linear maps whose adjoints are denoted as $\cA^*$ and $\cB^*$, respectively.
The vectors $l, u$ are given $p$-dimensional vectors whose elements are allowed to take the values $-\infty$ and $\infty$,
respectively.
Note that $\cP = \S^n$ is allowed, in which case there are no additional bound constraints imposed on $X$.
We assume that the $m\times m$ symmetric matrix $\cA\cA^*$ is invertible, i.e., $\cA$ is surjective.

Note that ({\bf SDP}) is equivalent to
\begin{eqnarray*}
({\rm P}) \quad  \min \Big\{ \inprod{C}{X} \mid \cA(X) = b,\; \cB(X) - s = 0, \; X \in \S^n_+, \; X\in \cP,  \; s \in \cQ \Big\},
\label{eq-primal}
\end{eqnarray*}
where $\cQ = \{s \in \Re^{p}: l \leq s \leq u\}$.
The dual of (P), ignoring the minus sign in front of the minimization, is given by
\begin{eqnarray*}
({\rm D})\; \min \left\{  \delta_{\cP}^*(-Z) + \delta_{\cQ}^*(-v)  + \inprod{-b}{y} \;\Big|\;
\begin{array}{l} \cA^*(y) + \cB^*(\bar{y}) + S + Z= C,\; -\bar{y} + v = 0,
\\[3pt]
 S\in \S^n_+,
\;Z\in \S^n, \; y\in \Re^m,\; \bar{y} \in \Re^p,\; v\in \Re^p
\end{array}
\right\},
\label{eq-dual}
\end{eqnarray*}
where for any $Z\in \S^n$, $ \delta_{\cP}^*(-Z) $ is defined by
\begin{eqnarray*}
\delta_{\cP}^*(-Z)
= \sup\{ \inprod{-Z}{W}\mid W \in \cP \}  \label{gz}
\end{eqnarray*}
and  $\delta_{\cQ}^*(\cdot)$ is defined similarly.
We note  that our solver is designed based on the assumption that (P) and (D) are feasible.

While we have presented the problem ({\bf SDP}) with a single variable block $X$, our
solver is capable of solving the following more general problem with $N$ blocks of  variables:
\begin{eqnarray}
 \begin{array}{ll}
 \min & \sum_{j=1}^N \inprod{C^{(j)}}{X^{(j)}}
\\[5pt]
\mbox{s.t.} & \sum_{j=1}^N \cA^{(j)}(X^{(j)}) = b, \quad
  l\leq \sum_{j=1}^N \cB^{(j)}(X^{(j)}) \leq u,
\\[5pt]
& X^{(j)}\in \cK^{(j)}, \; X^{(j)} \in \cP^{(j)}, \; j=1,\ldots,N,
\end{array}
\label{eq-gensdp}
\end{eqnarray}
where $\cA^{(j)}:\cX^{(j)} \to \Re^m$, and $\cB^{(j)}:\cX^{(j)} \to \Re^p$ are given linear maps,
$\cP^{(j)} := \{ X^{(j)}\in \cX^{(j)} \mid L^{(j)} \leq X^{(j)} \leq U^{(j)}\}$
and $L^{(j)}, U^{(j)}\in \cX^{(j)}$ are given symmetric matrices where the elements are allowed
to take the values $-\infty$ and $\infty$, respectively. Here
$\cX^{(j)} = \S^{n_j}$ ($\Re^{n_j}$), and $\cK^{(j)} = \cX^{(j)}$ or
$\cK^{(j)} = \S^{n_j}_+$  ($\Re^{n_j}_+$). For later expositions, we should note that when
$\cX^{(j)} = \S^{n_j}$,
the
linear map $\cA^{(j)}:\S^{n_j}\to \Re^m$ can be expressed in the form of
\begin{eqnarray}
  \cA^{(j)}(X^{(j)}) = \left[ \begin{array}{c}
  \inprod{A^{(j)}_1}{X^{(j)}}, \ldots,
  \inprod{A^{(j)}_m}{X^{(j)}}
\end{array} \right]^T,
\end{eqnarray}
where $A^{(j)}_1,\ldots,A^{(j)}_m \in \S^{n_j}$ are given constraint  matrices. The corresponding adjoint
$(\cA^{(j)})^*:\Re^m \to \S^{n_j}$ is then given by
\begin{eqnarray*}
\begin{array}{l}
  (\cA^{(j)})^* y = \sum_{k=1}^m y_k A^{(j)}_k.
\end{array}
\end{eqnarray*}

In this paper,
we introduce our {\sc Matlab} software
package \SDPNALP for solving ({\bf SDP}) or more generally \eqref{eq-gensdp},
where the maximum matrix dimension is assumed to be moderate (say less than 5000) but
the number of linear constraints $m+p$ can be large (say more than a million).
One of our main contributions here is that the current algorithm
has substantially extended the capability of \SDPNALP to solve
the general problem \eqref{eq-gensdp} compared to
the original version in \cite{sdpnalplus}, wherein the algorithm
is designed to solve a problem with only linear equality constraints
and $\cP = \{ X\in \S^n \mid X \geq 0\}$ or $\cP = \S^n$.
Moreover,
the implementation in  \cite{sdpnalplus} was based
on a majorized semismooth Newton-CG augmented Lagrangian method developed in
that paper. Here, for the purpose of deriving simpler stopping conditions, we
redesign  the algorithm by employing
an inexact  semi-proximal
alternating direction method of multipliers (sPADMM)
(or the semi-proximal augmented Lagrangian (sPALM) if the bound
constraints are absent)
framework
developed in \cite{inexactSPADMM}
for multi-block convex composite conic programming problems. Currently, the algorithm which we have implemented is
a $2$-phase algorithm  based on the
augmented Lagrangian function for (D).
In the first phase, we employ the inexact
symmetric Gauss-Seidel based sPADMM
to solve the problem to a modest level of accuracy.
Note that while
the main purpose of the first phase algorithm is to
generate a good initial point to warm-start
the  second phase algorithm, it can be
used on its own to solve a problem.
The algorithm we have implemented in the second phase is
an inexact sPADMM for which the main subproblem in each iteration
is solved by  a semismooth Newton-CG  method.

The development of \SDPNALP in \cite{sdpnalplus}, which is built on
the earlier work on {\sc Sdpnal} in \cite{ZST2010},
has in fact spurred
much of the recent progresses in designing efficient convergent ADMM-type
algorithms for solving multi-block convex composite conic programming,
such as \cite{inexactSPADMM,schurSPADMM,ADMM3c}.
Those works in turn shaped the recent algorithmic design of
\SDPNALPP
Indeed, the algorithm in the first phase of \SDPNALP is
the same as the convergent ADMM-type method developed in \cite{ADMM3c} when
the subproblems in each iteration are solved analytically.
For the algorithm in
the second phase, it is an economical variant of the
majorized semismooth Newton-CG (SNCG) augmented Lagrangian method
designed in \cite{sdpnalplus} to solve (D) for which
only one SNCG subproblem is solved in each iteration.

Another contribution of this paper is our development of a
basic interface for the users to input their SDP models into the \SDPNALP solver.
While there are currently two well developed {\sc matlab} based user interfaces
for SDP problems, namely,  CVX \cite{CVX} and YALMIP \cite{YALMIP}, there are strong motivations
for us to develop our own interface here. A new interface
 is necessary to
facilitate the modeling of an SDP problem for \SDPNALP  because
of latter's flexibility to directly accept inequality constraints of the form
``$l \leq {\cal B}(X) \leq u$'', and bound constraints of the form
``$L\leq X\leq U$''. The flexibility can significantly simplify the
generation of the data in the \SDPNALP format as compared
to what need to be done in  CVX or YALMIP
to reformulate them as equality constraints through introducing
extra variables.
In addition, the final number of  equality constraints
present in the data input to \SDPNALP can also
be substantially fewer than those present
in CVX or YALMIP. It is important to note here that the
number of equality constraints present in the generated problem data
can greatly affect the computational efficiency of the
solvers, especially for interior-point based solvers.
An illustration of the benefits just mentioned will be
given at the end of Section \ref{sec-interface-examples}.

Our \SDPNALP solver is designed for solving
feasible problems of the form presented in (P) and (D). It is capable
of solving large scale SDPs with $m$ or $p$ up to a few millions but $n$ is assumed
to be moderate (up to a few thousands). Extensive numerical experiments conducted in \cite{sdpnalplus} show that a variety of large scale SDPs can be solved by
\SDPNALP much more efficiently than the best alternative methods
\cite{monteiro2013first,wen2010alternating}.

The \SDPNALP package can be downloaded from the following website:
\begin{center}
\url{http://www.math.nus.edu.sg/\~mattohkc/SDPNALplus.html}
\end{center}
Installation and general information such as citations, can
be found at the above link.
The test instances which we have used to evaluate the performance of
our solver can also be found at the above website.

We have evaluated the performance of \SDPNALP on various classes of large scale SDP
problems arising from the relaxation of combinatorial problems such as maximum stable
set problems, quadratic assignment problems, frequency assignment problems, and
binary integer quadratic programming problems. The solver has also been tested
on large SDP problems arising from robust clustering problems, rank-one tensor
approximation problems, as well as electronic structure calculations in quantum
chemistry. The detailed numerical results can be found at the above website.
Based on the numerical evaluation of \SDPNALP on 745  SDP problems,
we can observe that the solver is fairly robust (in the sense that
it is able to solve most of the tested problems to the accuracy of
$10^{-6}$ in the relative KKT residual) and highly efficient
in solving the tested classes of problems.

The remaining parts of this paper are organized as follows.
In the next section, we describe the installation and present some general
information on our software.
Section \ref{sec-imple} gives some details on
the main solver function {\tt sdpnalplus.m}.
In Section \ref{sec-algo}, we describe the algorithm implemented in
\SDPNALP and discuss some implementation issues.
In Section \ref{sec-interface}, we present a
basic interface for the users to input their SDP models into the \SDPNALP solver.
In Section \ref{sec-interface-examples},

we present a few SDP examples to illustrate the usage of our software,
and how to
input the SDP models into our interface.
Section \ref{sec-num} gives a summary of the numerical
results obtained by \SDPNALP in solving 745 test instances of SDP problems
arising from various sources.
Finally, we conclude
the paper in Section \ref{sec-conclusion}.

\
%
%
%

\section{Data structure and main solver} \label{sec-imple}

\SDPNALP is an enhanced version of the {\sc Sdpnal} solver
developed by Zhao, Sun and Toh \cite{ZST2010}.
The internal implementation of  \SDPNALP thus follows the data structures
and design framework of {\sc Sdpnal}. A casual user need not
understand the internal implementation of \SDPNALP\!\!.

\subsection{The main function: {\tt sdpnalplus.m}}

In the \SDPNALP solver, the main routine is {\tt sdpnalplus.m}, whose calling syntax is as follows:
\begin{verbatim}
[obj,X,s,y,S,Z,ybar,v,info,runhist] = ...
     sdpnalplus(blk,At,C,b,L,U,Bt,l,u,OPTIONS,X,s,y,S,Z,ybar,v);
\end{verbatim}
{\bf Input arguments.}

\begin{itemize}
\item {\tt blk}: a cell array describing the conic block structure of the SDP problem.

\item {\tt At, C, b, L, U, Bt, l, u}: data of the problem ({\bf SDP}).
\\
If ${L\leq X}$ but $X$ is unbounded above, one can set {\tt U=inf} or {\tt U=[]}.
Similarly, if the linear map $\cB$ is not present, one can set {\tt Bt=[]}, {\tt l=[]},
{\tt u=[]}.

\item {\tt OPTIONS}: a structure array of parameters (optional).

\item {\tt X, s, y, S, Z, ybar, v}: an initial iterate (optional).

\end{itemize}

\noindent {\bf Output arguments.}
The names chosen for the output arguments explain their contents. The
argument ${\tt X}$ is a solution to (P) which satisfies the constraints
$X \in \S^n_+$ and
 $X\in\cP$  approximately up to the desired accuracy tolerance.
The argument {\tt info} is a structure array which records various performance measures
of the solver. For example
\begin{verbatim}
  info.etaRp, info.etaRd, info.etaK1, info.etaK2
\end{verbatim}
correspond to the measures $\eta_P$, $\eta_D$, $\eta_{\cK}$, $\eta_{\cP}$ defined
later
in \eqref{eq-eta}, respectively.
The argument {\tt runhist} is a structure array which records the history of various performance
measures during
the course of running {\tt sdpnalplus.m}. For example,
\begin{eqnarray*}
& {\tt runhist.primobj},\;\; {\tt runhist.dualobj},\;\; {\tt runhist.relgap} &
\\
& {\tt runhist.primfeasorg},\;\; {\tt runhist.dualfeasorg} &
\end{eqnarray*}
record the primal and dual objective values, complementarity gap, primal and dual infeasibilities at each iteration, respectively.

\subsection{Generation of starting point by {\tt admmplus.m} }

If an initial point {\tt (X,s,y,S,Z,ybar,v)} is not provided
for {\tt sdpnalplus.m}, we call
the function {\tt admmplus.m}, which implements a convergent $3$-block ADMM proposed in \cite{ADMM3c},  to generate a starting point.
The routine {\tt admmplus.m} has a similar calling syntax as {\tt sdpnalplus.m} given as follows:
\begin{verbatim}
[obj,X,s,y,S,Z,ybar,v,info,runhist] = ...
     admmplus(blk,At,C,b,L,U,Bt,l,u,OPTIONS,X,s,y,S,Z,ybar,v);
\end{verbatim}
Note that if an initial point {\tt (X,s,y,S,Z,ybar,v)} is not supplied to {\tt admmplus.m},
the default initial point is {\tt (0,0,0,0,0,0,0)}.

We should mention that although we use {\tt admmplus.m}  for the
purpose of warm-starting {\tt sdpnalplus.m}, the user has the freedom to use {\tt admmplus.m}
alone to solve the problem ({\rm \bf SDP}).

\subsection{Arrays of input data}
\label{subsec-data}

The format of the input data in \SDPNALP is similar to those in SDPT3 \cite{SDPT3a,SDPT3b}. For each SDP problem,  the conic block structure of the problem data is described by a cell array named {\tt blk}. If the $k$th block {\tt X\{k\}} of the variable {\tt X}  is a nonnegative vector block with dimension $n_k$, then we set
\begin{enumerate}
\item[] {\tt blk\{k,1\} = 'l'},\; {\tt blk\{k,2\} =} $n_k$,

\item[] {\tt At\{k\} =} [$n_k \times m$ \, sparse], \;
{\tt Bt\{k\} =} [$n_k \times p$ \, sparse],

\item[] {\tt C\{k\}, L\{k\}, U\{k\}, X\{k\}, S\{k\}, Z\{k\} =}
 [$n_k \times 1$ \, double or sparse].
\end{enumerate}
If the $j$th block {\tt X\{j\}} of the variable {\tt X}
 is a semidefinite block consisting of a single block of size $s_j$, then the content of the $j$th block is given as follows:
\begin{enumerate}
\item[] {\tt blk\{j,1\} = 's'}, \; {\tt blk\{j,2\} =} $s_j$,

\item[] {\tt At\{j\} =} [$\bar{s}_j \times m $ \,  sparse ],\;
{\tt  Bt\{k\} =} [$\bar{s}_j \times p$ \, sparse],

\item[] {\tt C\{j\}, L\{j\}, U\{j\}, X\{j\}, S\{j\}, Z\{j\} =}
 [$s_j\times s_j$  \, double or sparse ],
\end{enumerate}
where $\bar{s}_j = s_j(s_j+1)/2$.
By default, the contents of the cell arrays {\tt L} and {\tt U} are set to be empty arrays.
But if ${\tt X\{j\}} \geq 0$ is required, then one
can set
$$
 \mbox{\tt L\{j\} = 0}, \quad \mbox{\tt U\{j\} = []}.
$$
One can also set $L=0$ to indicate that ${\tt X\{j\}} \geq 0$ for all $j=1,\ldots,N$
in \eqref{eq-gensdp}.

We should mention that for the sake of computational efficiency, we store all the constraint matrices associated with the $j$th semidefinite block in vectorized form as a single $\bar{s}_j \times m $
matrix {\tt At\{j\}}, where the $k$th column of this matrix corresponds to the $k$th constraint matrix $A^{(j)}_k$, i.e.,
\[
{\tt At\{j\}} = [{\rm svec}(A^{(j)}_{1}), \ldots, {\rm svec} (A^{(j)}_{m})],
\]
and ${\rm svec}: \mathcal{S}^{s_j} \rightarrow \Re^{\bar{s}_j}$ is the vectorization operator on symmetric matrices defined by
\begin{eqnarray}\label{fun-svec}
{\rm svec}(X) = [X_{11},\sqrt{2}X_{12},X_{22}, \ldots, \sqrt{2}X_{1,s_j},\ldots, \sqrt{2}X_{s_j-1,s_j}, X_{s_j,s_j}]^T .
\end{eqnarray}
We store {\tt Bt} in the same format as {\tt At}. The function {\tt svec.m} provided in \SDPNALP can easily convert a symmetric matrix into the vector storage scheme
described in \eqref{fun-svec}.
Note that while we store the constraint matrices in vectorized form,
the semidefinite blocks in the
variables {\tt X}, {\tt S} and {\tt Z}
are stored either as matrices or in vectorized forms according to the storage scheme of the
input data {\tt C}.

Other than inputting the data {\tt (At,b,C,L,U)} of an SDP problem individually,
\SDPNALP also provides the
functions {\tt read\_sdpa.m} and {\tt read\_sedumi.m} to convert
problem data stored in the
SDPA \cite{SDPA} and SeDuMi \cite{SeDuMi} format   into our cell-array data format just described. For example, for the problem {\tt theta62.dat-s}
in the folder {\tt /datafiles}, the user can call the m-file {\tt read\_sdpa.m} to load
the SDP data as follows:
\begin{small}
\begin{verbatim}
>> [blk,At,C,b] = read_sdpa('./datafiles/theta62.dat-s');
>> OPTIONS.tol = 1e-6;
>> [obj,X,s,y,S,Z,ybar,v,info,runhist] = sdpnalplus(blk,At,C,b,[],[],[],[],[],OPTIONS);
\end{verbatim}
\end{small}

\subsection{The structure array {\tt OPTIONS} for parameters}
\label{subsec-para}
Various parameters used in our solver {\tt sdpnalplus.m} are set
in the structure array {\tt OPTIONS}. For details, see
{\tt SDPNALplus\_parameters.m}.
The important parameters which the user is likely to reset are described next.
\begin{enumerate}
\item {\tt OPTIONS.tol}: accuracy tolerance to terminate the algorithm, default is $10^{-6}$.

\item {\tt OPTIONS.maxiter}: maximum number of iterations allowed, default is $20000$.

\item {\tt OPTIONS.maxtime}: maximum time (in seconds) allowed, default is $10000$.

\item {\tt OPTIONS.tolADM}: accuracy tolerance to use for {\tt admmplus.m} when                            generating a starting point  for the algorithm in the second phase
of {\tt sdpnalplus.m} (default = $10^{-4}$).

\item {\tt OPTIONS.maxiterADM}: maximum number of ADMM iterations allowed for
generating a starting point.
When there are no bound constraints on $X$ ($\cP = \S^n$) and no linear inequality constraints
corresponding to $\cB(X)$ (hence $\cQ = \emptyset$),
the default value is roughly equal to 200; otherwise, the default value is 2000.

\item {\tt OPTIONS.printlevel}: different levels of details to print the intermediate information during the run.
It can be the integers $0,1,2$, with $1$ being the default.
Setting to the highest value $2$ will result in printing the complete details.

\item {\tt OPTIONS.stopoption}: options to stop the solver.
The default is {\tt OPTIONS.stopoption=1}, for which the solver may be stopped prematurely
when stagnation occurs. To prevent the solver from stopping prematurely before
the required accuracy is attained, set {\tt OPTONS.stopoption=0}.

\item{ \tt OPTIONS.AATsolve.method}: options to solve a
linear system involving the coefficient matrix $\cA\cA^*$, with
\\[5pt]
{\tt
OPTIONS.AATsolve.method='direct' (default) or  'iterative'. }
\\[5pt]
For the former option, a linear system of the form $\cA\cA^* y=h$ is solved by the sparse Cholesky factorization, while for the latter option, it is solved by a diagonally preconditioned PSQMR
iterative solver.
\end{enumerate}

\subsection{Stopping criteria}

In \SDPNALP\!\!,   we measure the accuracy of an approximate optimal solution $(X,s,y,\by,S,Z,v)$ for (P) and (D) by using the following relative residual based on the KKT optimality conditions:
\begin{eqnarray}
\eta = \max \{\eta_P,\eta_D,\eta_{\cK},\eta_{\cP}\},
 \label{eq-eta}
\end{eqnarray}
where $\cK = \S^n_+$,
\begin{eqnarray*}
\begin{array}{l}
\eta_P = \max\Big\{\frac{\norm{\cA(X)-b}}{1+\norm{b}}, \frac{\norm{\cB(X)-s}}{1+\norm{s}}\Big\},\;
\eta_D =\max\Big\{\frac{\norm{\cA^*(y)+\cB^*(\bar{y})+S +Z-C}}{1+\norm{C}},\frac{\norm{\bar{y}-v}}{1+\norm{v}}\Big\},
\\[8pt]
\eta_{\cK} = \frac{1}{5} \frac{\norm{X-\Pi_{\cK}(X-S)}}{1+\norm{X}+\norm{S}},\;
\eta_{\cP} = \frac{1}{5}\max\Big\{\frac{\norm{X-\Pi_{\cP}(X-Z)}}{1+\norm{X}+\norm{Z}},\frac{\norm{s-\Pi_{\cQ}(s-v)}}{1+\norm{s}+\norm{v}}\Big\}.
\end{array}
\end{eqnarray*}
Additionally, we compute the relative gap by
\begin{eqnarray}
&\eta_g = \frac{|{\tt pobj} - {\tt dobj}|}{1+\abs{\tt pobj}+\abs{\tt dobj}}.&
\label{eq-gap}
\end{eqnarray}
For a given  accuracy tolerance specified in  {\tt OPTIONS.tol},
we terminate both {\tt sdpnalplus.m} and {\tt admmplus.m} when
\begin{eqnarray}
\eta  \leq {\tt OPTIONS.tol}. \label{stop}
\end{eqnarray}

\subsection {Caveats}
There are a few points which we should emphasize on our solver.
\begin{itemize}
\item It is important to note that \SDPNALP is a research software. It is not intended
nor designed to be a general purpose software at the moment.
The solver is designed based on the assumption that the primal and dual
SDP problems (P) and (D) are feasible, and that Slater's constraint qualification holds.
 The solver is expected to be robust if the primal and dual SDP problems
are both non-degenerate at the optimal solutions. However, if either
 one of them, particularly if the primal problem,  is degenerate or if the Slater's condition fails, then the solver may not be able to solve
the problems to high accuracy.

\item Another point to note is that our solver is designed with the emphasis on
handling problems with positive semidefinite variables efficiently. Little attention has been
paid on optimizing the solver to handle linear programming problems.

\item While in theory our solver can easily be extended to solve problems with
second-order cone constraints,  it is not capable of solving such problems at the moment
although we plan to extend our solver to handle second-order cone programming
problems in the future.

\end{itemize}

\section{Algorithmic design and implementation}
\label{sec-algo}

For simplicity, we will describe the algorithmic design for the
problem (D) instead of the dual of the more general problem \eqref{eq-gensdp}.
Our algorithm is developed based on the augmented Lagrangian function
for (D), which is defined as follows: given a penalty parameter $\sig > 0$,
for $(Z,v,y,\bar{y})\in \S^n\times \Re^p\times \times
\Re^m\times \Re^p$, and $(X,s)\in\S^n\times \Re^p$,
\begin{eqnarray*}
 L_\sig(Z,v,S,y,\bar{y}; X,s) = \left\{ \begin{array}{l}
\delta_{\cP}^*(-Z) + \delta_{\cQ}^*(-v)  + \inprod{-b}{y} + \delta_{\S^n_+}(S)
-\frac{1}{2\sig}\norm{X}^2 -\frac{1}{2\sig}\norm{s}^2
\\[5pt]
+\frac{\sig}{2}\norm{\cA^*(y) + \cB^*(\bar{y}) + S + Z-C+\sig^{-1} X}^2
+\frac{\sig}{2}\norm{v-\bar{y} +\sig^{-1} s}^2 .
\end{array}\right.
\end{eqnarray*}

As mentioned in the Introduction, the algorithm implemented in
\SDPNALP is a 2-phase algorithm where the first phase is a convergent inexact sGS-sPADMM algorithm \cite{inexactSPADMM}
whose  template is described next.

\bigskip
\noindent{\bf First-phase algorithm.}
Given an initial iteration $(Z^0,v^0,S^0,y^0,\by^0,X^0,s^0)$,
perform the following steps in each iteration.
\begin{description}
\item[Step 1.]  Let $R_{1}^k = \cA^*(y^k) + \cB^*(\by^k) + S^k + Z^k-C+\sig^{-1} X^k$
and $R^k_2 = v^k- \by^k+\sig^{-1}s^k$.
Compute $(Z^{k+1},v^{k+1}) = \mbox{argmin}\;  L_\sig(Z,v,S^k,y^{k},\by^{k};X^k,s^k)$ as follows:
\begin{eqnarray*}
  Z^{k+1} &=&  \mbox{argmin}\big\{ \delta_{\cP}^*(-Z)
+\frac{\sig}{2}\norm{Z-Z^k+R_1^k}^2   \big\} =  \sig^{-1} \Pi_{\cP}(\sig(R^k_1-Z^k)) - (R^k_1-Z^k ),
\\[5pt]
 v^{k+1} &=&  \mbox{argmin}\big\{  \delta_{\cQ}^*(-v)
+\frac{\sig}{2}\norm{v-v^k+R_2^k}^2   \big\} = \sig^{-1} \Pi_{\cQ}(\sig(R^k_2-v^k)) - (R^k_2-v^k ).
\end{eqnarray*}

\item[Step 2a.]
Compute
\begin{eqnarray*}
 (y^{k+1}_{\rm tmp},\by^{k+1}_{\rm tmp}) & \approx &\mbox{argmin} \big\{ L_\sig(Z^{k+1},v^{k+1},S^k,y,\by;X^k,s^k)
\big\}.
\end{eqnarray*}
 For this step,
we typically need to solve a large system of linear equations given by
\begin{equation}
 \underbrace{\left[\begin{array}{cc}
 \cA\cA^* & \cA\cB^* \\[5pt] \cB\cA^* & \cB\cB^*+\cI
\end{array}\right]}_{\cM}
\left[\begin{array}{c} y \\[5pt]\by \end{array}\right]
= \left[\begin{array}{l} h_1 := \sig^{-1}b -{\cA(S^k+Z^{k+1}-C+\sig^{-1}X^k)}
\\[5pt] h_2:=v^{k+1}+\sig^{-1}s^k -\cB(S^k+Z^{k+1}-C+\sig^{-1}X^k)
\end{array}\right].
\label{eq-step2a}
\end{equation} 
In our implementation, we solve the linear system via the sparse Cholesky factorization
of $\cM$ if it  can be computed at a moderate cost. Otherwise,
we use a preconditioned CG method to solve \eqref{eq-step2a}
approximately so that the residual norm satisfies the following accuracy condition:
\begin{eqnarray*}
 \sqrt{\sig} \norm{[h_1; h_2] - \cM [y^{k+1}_{\rm tmp}; \by^{k+1}_{\rm tmp}] }\leq
 \varepsilon_k,
\end{eqnarray*}
where $\{ \varepsilon_k\}$ is a predefined summable sequence of nonnegative numbers.
In \cite{sdpnalplus}, the linear system corresponding to $\cM$ is $\cA\cA^* y = h_1$,
and it is solved by a direct method based on sparse Cholesky factorization.
Here, the inexact sGS-ADMM framework \cite{inexactSPADMM} we have employed
gives us the flexibility to solve the linear system  approximately by an
iterative solver such as the preconditioned conjugate gradient method, while not affecting the convergence of the algorithm.
Such a flexibility is obviously critical to the computational efficiency of the
algorithm when the sparse Cholesky
factorization of $\cM$ is impossible to compute for a very large linear system.
\item[Step 2b.] Let $R^{k+1}_1 = \cA^* (y^{k+1}_{\rm tmp})
+ \cB^* (\by^{k+1}_{\rm tmp})
+ S^k + Z^{k+1} - C + \sig^{-1} X^k$.
Compute
\begin{eqnarray*}
  S^{k+1} &=&
 \mbox{argmin} \big\{ \delta_{\S^n_+}(S)+\frac{\sig}{2}\norm{S - S^k + R_1^{k+1}} ^2
\big\}
= \Pi_{\S^n_+} ( S^k-R_1^{k+1}).
\end{eqnarray*}
\item[Step 2c.]
Let $h_1^{\rm new} := h_1 - \cA(S^{k+1}-S^k)$, and $h_2^{\rm new} := h_2 - \cB(S^{k+1}-S^k)$.
Set $(y^{k+1},\by^{k+1}) = (y^{k+1}_{\rm tmp},\by^{k+1}_{\rm tmp})$
if
\begin{eqnarray*}
  \sqrt{\sig} \left\| [h_1^{\rm new}; {h_2^{\rm new}}] -\cM [y^{k+1}_{\rm tmp};
\by^{k+1}_{\rm tmp}] \right\| \leq 10  \varepsilon_k;
\end{eqnarray*}
otherwise solve \eqref{eq-step2a} with the vector $h_1$ replaced
by $h_1^{\rm new}$ and $h_2$ replaced by $h_2^{\rm new}$, and
the approximate solution $(y^{k+1},\by^{k+1})$ should satisfy the above
accuracy condition.
\item[Step 3.] Let
$R_{D,1}^{k+1} = \cA^*(y^{k+1}) + \cB^*(\by^{k+1}) + S^{k+1} + Z^{k+1}-C$
and $R^{k+1}_{D,2} = v^{k+1}- \by^{k+1}$.
Compute
$$
 X^{k+1} = X^k + \tau \sig R^{k+1}_{D,1}, \; s^{k+1} = s^k + \tau\sig R^{k+1}_{D,2},
$$
where $\tau \in (0,(1+\sqrt{5})/2)$ is the steplength which is typically
chosen to be 1.618.
\end{description}

\bigskip
We note that by \cite{inexactSPADMM}, the computation in Step 2a--2c is equivalent to solving
the subproblem:
\begin{eqnarray*}
(S^{k+1},y^{k+1},\bar{y}^{k+1}) =  \mbox{argmin} \left\{
\begin{array}{l}
{ L_\sig(Z^{k+1},v^{k+1},S,y,\bar{y}; X^k,s^k)}
\\[5pt]
 + \frac{\sig}{2}\norm{(S;y;\bar{y}) - (S^k;y^k;\bar{y}^k)}^2_{\cH}
 \end{array}
 \right\},
\end{eqnarray*} 
where $\cH$ is the symmetric Gauss-Seidel decomposition linear operator associated with
the linear operator $(\cI;\cA;\cB) (\cI,\cA^*,\cB^*) + \diag(0,0,\cI)$, i.e.,
\begin{eqnarray*}
  \cH = \left[  \begin{array}{ccc}
   (\cA^*,\, \cB^*) \cD^{-1} (\cA; \cB) & 0 & 0 \\[5pt]
    0 & 0 & 0 \\[5pt]
    0 & 0 & 0
  \end{array}
  \right] \quad { \mbox{with}} \;\; \cD = \left[\begin{array}{cc}\cA\cA^* & \cA\cB^*\\ \cB\cA^* & \cB\cB^*+\cI\end{array} \right].
\end{eqnarray*}

\bigskip
There are numerous implementation issues which are addressed in
\SDPNALP to make the above skeletal algorithm practically efficient and robust.
A detailed description of how the issues are addressed is beyond the scope of this paper.
Hence we shall only briefly mention the most crucial ones.
\begin{enumerate}
\item Dynamic adjustment of the penalty parameter $\sig$, which is equivalent
to restarting the algorithm with a new parameter by using the
most recent iterate as the initial starting point.
\item Initial scaling of the data, and dynamic scaling of the data.
\item The efficient  implementation of the PCG method to compute an approximate solution for \eqref{eq-step2a}.
\item Efficient computation of the iterate $S^{k+1}$ by using partial eigenvalue decomposition
whenever it is expected to be more economical than a full eigenvalue decomposition.
\item Efficient evaluation of the residual measure $\eta$ defined in \eqref{eq-eta}.
\end{enumerate}

The algorithm in the second phase of \SDPNALP is designed based on the following
convergent
inexact sPADMM algorithm (or the sPALM algorithm if the bound constraints are absent).
After presenting the algorithm, we will explain the
changes we made in this algorithm compared to that developed
in \cite{sdpnalplus}.

\bigskip
\noindent{\bf Second-phase algorithm.}
Given an initial iterate $(Z^0,v^0,S^0,y^0,\by^0,X^0,s^0)$
 generated in the first phase, perform the following steps in each iteration.
\begin{description}
\item[Step 1.]  Compute $(Z^{k+1},v^{k+1})$ as in Step 1 of the first-phase algorithm.
\item[Step 2.] Compute
\begin{eqnarray*}
 (y^{k+1},\by^{k+1},S^{k+1}) \approx \mbox{argmin} \, L_\sig (Z^{k+1},v^{k+1},S,y,\by;X^k,s^k)
\end{eqnarray*}
by using the semismooth Newton-CG (SNCG) method which has been described in detail in \cite{ZST2010}
such that the following accuracy condition is met:
\begin{eqnarray*}
  \sqrt{\sig}\max\{\norm{ b - \cA \Pi_{\S^n_+}(W^{k+1})},
   \norm{ \cB \Pi_{\S^n_+}(W^{k+1})-s^k+\sig(\by^{k+1}-v^{k+1})} \} \leq \varepsilon_k,
\end{eqnarray*}
where
$W^{k+1} := \cA^* y^{k+1}+\cB^*\by^{k+1}+S^k+Z^{k+1}-C+\sig^{-1}X^k$,
and $\{ \varepsilon_k\}$ is a predefined summable sequence of nonnegative numbers.
\item[Step 3.]
Compute $(X^{k+1},s^{k+1})$ as in Step 3 of the first-phase algorithm.
\end{description}

\bigskip
As one may observe, the difference between the first-phase and the second-phase
algorithms lies in the construction of
$(y^{k+1},\by^{k+1},S^{k+1})$ in Step 2 of the algorithms.
In the first phase, the iterate is generated by adding the
semi-proximal term
$
 \frac{\sig}{2}\norm{(S;y;\bar{y})-(S^k;y^k;\bar{y}^k)}^2_{\cH}$
 to the
augmented Lagrangian function $L_\sig(Z^{k+1},v^{k+1},S,y,\by;X^k,s^k)$.
For the second phase, no such a semi-proximal term is required
though one may still add a small semi-proximal term
to the augmented Lagrangian function to ensure that the subproblems
are well defined.
As our goal is to minimize the augmented Lagrangian function $L_\sig(Z,v,S,y,\bar{y};X^k,s^k)$
for each pair of given $(X^k,s^k)$, it is thus clear that Step 2 of the second-phase algorithm
is closer to that goal compared to Step 2 of the first-phase algorithm.
Of course, the price to pay is that the subproblem in Step 2 of the second-phase algorithm
is more complicated to solve.

Now we highlight the differences between the above inexact sPADMM algorithm
and the majorized semismooth Newton-CG (MSNCG) augmented Lagrangian method  developed in \cite{sdpnalplus}. First, the algorithm in \cite{sdpnalplus} is designed to
solve ({\bf SDP}) with only linear equality constraints while the algorithm here
is for the general problem with additional linear inequality constraints.
Even when we specialize the algorithm here to the problem with only
linear equality constraints, our algorithm here is also different from the one in \cite{sdpnalplus}
which we will now explain. For the case when only linear equality constraints are present,
the augmented Lagrangian function associated with the dual of that problem is given by
\begin{eqnarray*}
 L_\sig(Z,S,y; X) \;=\; \delta_\cP^*(-Z) + \inprod{-b}{y} + \delta_{\S^n_+}(S)
 +\frac{\sig}{2} \norm{\cA^* y + S + Z - C +\sig^{-1} X}^2 -\frac{1}{2\sig}\norm{X}^2.
\end{eqnarray*}
At the $k$th iteration of the MSNCG augmented Lagrangian method, the following
subproblem must be solved:
\begin{eqnarray*}
 \min_{y,S,Z} \big\{ L_\sig(Z,S,y; X^k) \big\},
\end{eqnarray*}
and theoretically it is solved by the MSNCG method until a certain stopping condition is
satisfied. However, in the practical implementation, only one step of the MSNCG method
is applied to solve the subproblem and the stopping condition is not {strictly enforced}. Thus there is a gap between
the theoretical algorithm and the practical algorithm implemented in
\cite{sdpnalplus}. But for the convergent inexact sPADMM algorithm employed in this paper,
its practical implementation  follows closely the steps described
in the second-phase algorithm. Thus the practical algorithm presented
in this paper is based on rigorous stopping conditions in each iteration to guarantee its overall convergence.


\section{Interface}
\label{sec-interface}

In this section, we will present a basic interface for our \SDPNALP solver. First, we show how to use it via a small SDP example given as follows:
\begin{equation}\label{EX1-interface}
\begin{array}{rl}
\min & {\rm trace}(X^{(1)}) + {\rm trace}(X^{(2)}) + {\rm sum}(X^{(3)})\\[5pt]
\rm{s.t.} & -X^{(1)}_{12} + 2X^{(2)}_{33} + 2X^{(3)}_{2} = 4,\\[5pt]
& 2X^{(1)}_{23} + X^{(2)}_{42} - X^{(3)}_{4} = 3,\\[5pt]
& 2 \leq -X^{(1)}_{12} - 2X^{(2)}_{33} + 2X^{(3)}_{2} \leq 7, \\ [5pt]
& X^{(1)} \in \S^{6}_{+}, \; X^{(2)} \in \Re^{5\times 5}, \; X^{(3)} \in \Re^{7}_{+}, \\ [5pt]
& 0\leq X^{(1)} \leq 10E_{6}, \; 0\leq X^{(2)} \leq 8E_{5},
\end{array}
\end{equation}
where $E_n$ denotes the $n\times n$ matrix of all ones. In the notation of \eqref{eq-gensdp}, the problem \eqref{EX1-interface} has three blocks of variables $X^{(1)}$, $X^{(2)}$, $X^{(3)}$. The first linear map $\mathcal{A}^{(1)}$ contains two constraint matrices $A^{(1)}_1, A^{(1)}_2 \in \S^{6}$ whose nonzero elements are given by
\begin{equation*}
\begin{array}{cc}
(A^{(1)}_1)_{12} = (A^{(1)}_1)_{21} = -0.5, & (A^{(1)}_2)_{23} = (A^{(1)}_2)_{32} = 1.
\end{array}
\end{equation*}
With the above constraint matrices, we get $\langle A^{(1)}_1, X^{(1)}\rangle = -X^{(1)}_{12}$ and $\langle A^{(1)}_2, X^{(1)}\rangle = 2X^{(1)}_{23}$. \\
The second linear map $\mathcal{A}^{(2)}$ contains two constraint matrices $A^{(2)}_1, A^{(2)}_2 \in \Re^{5\times 5}$ whose nonzero elements are given by
\begin{equation*}
\begin{array}{cc}
(A^{(2)}_1)_{33} = 2, & (A^{(2)}_2)_{42} = 1.
\end{array}
\end{equation*}
Since the third variable $X^{(3)}$ is a vector, the third linear map $\mathcal{A}^{(3)}$ is a {constraint} matrix $A^{(3)} \in \Re^{2\times 7}$ whose nonzero elements are given by
\begin{equation*}
\begin{array}{cc}
(A^{(3)})_{12} = 2, & (A^{(3)})_{24} = -1.
\end{array}
\end{equation*}
In a similar fashion, one can identify the matrices for the linear maps $\mathcal{B}^{(1)}, \mathcal{B}^{(2)}$, and $\mathcal{B}^{(3)}$.

The example \eqref{EX1-interface} can be coded using our interface as follows:
\begin{footnotesize}
\begin{tt}
\begin{lstlisting}[caption = {Example \eqref{EX1-interface}.}\label{code: list1}, frame = single]
n1 = 6; n2 = 5; n3 = 7;
mymodel = ccp_model('Example_simple');
X1 = var_sdp(n1,n1);
X2 = var_nn(n2,n2);
X3 = var_nn(n3);
mymodel.add_variable(X1,X2,X3);
mymodel.minimize(trace(X1) + trace(X2) + sum(X3));
mymodel.add_affine_constraint(-X1(1,2)+2*X2(3,3)+2*X3(2) == 4);
mymodel.add_affine_constraint(2*X1(2,3)+X2(4,2)-X3(4) == 3);
mymodel.add_affine_constraint(2<=-X1(1,2)-2*X2(3,3)+2*X3(2)<=7);
mymodel.add_affine_constraint(0 <= X1 <= 10);
mymodel.add_affine_constraint(X2 <= 8);
mymodel.solve;
\end{lstlisting}
\end{tt}
\end{footnotesize}
Note that although the commands
\begin{footnotesize}
\begin{verbatim}
mymodel.add_affine_constraint(-X1(1,2)+2*X2(3,3)+2*X3(2)==4);
mymodel.add_affine_constraint(2*X1(2,3)+X2(4,2)-X3(4)==3);
\end{verbatim}
\end{footnotesize}
are convenient to use for a small example, it may become tedious if there are many such constraints. In general, it is more economical to encode numerous such constraints by using the constraint matrices of the linear maps
$\mathcal{A}^{(1)}$, $\mathcal{A}^{(2)}$, $\mathcal{A}^{(3)}$, which we illustrate below:
\begin{footnotesize}
\begin{tt}
\begin{lstlisting}[caption = {Example  \eqref{EX1-interface} with constraints specified via linear maps as cell arrays.}\label{code: init_as_cell},frame=single]
A1 = {sparse(n1,n1); sparse(n1,n1)}; A2 = {sparse(n2,n2); sparse(n2,n2)};
A3 = sparse(2,n3);
A1{1}(1,2) = -1; A2{1}(3,3) = 2; A3(1,2) = 2;  % -X1(1,2)+2*X2(3,3)+2*X3(2)
A1{2}(2,3) = 2;  A2{2}(4,2) = 1; A3(2,4) = -1; % 2*X1(2,3)+X2(4,2)-X3(4)
b = [4;3];
mymodel.add_affine_constraint(A1*X1 + A2*X2 + A3*X3 == b);
\end{lstlisting}
\end{tt}
\end{footnotesize}
As the reader may have noticed, in constructing the matrix \texttt{A1\{1\}} corresponding to the constraint
matrix $A^{(1)}_1$, we set {\tt A1\{1\}(1,2) = -1} instead of {\tt A1\{1\}(1,2) = -0.5; A1\{1\}(2,1) = -0.5}.
Both ways of inputing {\tt A1\{1\}} are acceptable as internally, we will symmetrize the matrix {\tt A1\{1\}}.

In following subsections, we will discuss the details of the interface.

\subsection{Creating a ccp model}

Before declaring \texttt{variables}, \texttt{constraints} and setting \texttt{parameters}, we need to
create a {\tt ccp\_model} class  first.
This is done via the command:
\begin{center}
{\texttt{mymodel}} = \textbf{ccp\_model}(\texttt{model\_name});
\end{center}
The string \texttt{model\_name} is the name of the created {\tt ccp\_model}. If no model name is
specified, the default name is \texttt{`Default'}.

After solving the created \texttt{mymodel}, we save all the relevant information in the file
  \texttt{`model\_name.mat'}. It contains two structure arrays, \texttt{input\_data} and  \texttt{solution},
  which
 store all the input data and solution information, respectively.

\subsection{Delcaring variables}
\label{subsec-var}

Variables in \SDPNALP can be real vectors or matrices. Currently, our interface supports four types of {variables}: {free variables}, {variables in SDP cones}, {nonnegative variables} and
{variables which are symmetric matrices}. Next, we introduce them in details.

\begin{description}
\item[1.] \textbf{Free variables.} One can declare a  free variable ${\tt X} \in \Re^{m\times n}$  via the command:
\begin{eqnarray*}
 \mbox{{\tt X} = \textbf{var\_free}({\tt m},{\tt n}); }
\end{eqnarray*}
where the parameters {\tt m} and {\tt n} specify the dimensions of {\tt X}.
One can also declare a column vector variable ${\tt Y}\in\Re^n$ simply via the command:
\begin{eqnarray*}
 \mbox{{\tt Y} = \textbf{var\_free}({\tt n}); }
\end{eqnarray*}
\item[2.] \textbf{Variables in  SDP cones.}
A variable ${\tt X} \in \S^n_+$ can be declared via the command:
\begin{eqnarray*}
\mbox{{\tt X} = \textbf{var\_sdp}({\tt n},{\tt n});}
\end{eqnarray*}
In this case, the variable must be a square matrix,
so {\tt X} = \textbf{var\_sdp}({\tt m},{\tt n}) with ${\tt m} \not = {\tt n}$ is invalid.
\item[3.] \textbf{Variables in nonnegative orthants.}  To declare a nonnegative variable
${\tt X} \in \Re^{m\times n}_+$, one can use the command:
\begin{eqnarray*}
\mbox{{\tt X} = \textbf{var\_nn}({\tt m},{\tt n});}
\end{eqnarray*}
We can also use {\tt Y} = \textbf{var\_nn}({\tt n}) to declare a
vector variable ${\tt Y} \in \Re^n_+$.
\item[4.] \textbf{Variables which are symmetric matrices.}
 To declare a symmetric matrix variable ${\tt X} \in \S^n$, one can use the command:
\begin{eqnarray*}
 \mbox{{\tt X} = \textbf{var\_symm}({\tt n},{\tt n});}
\end{eqnarray*}
In this case, the variable must be a square matrix.

\item[5.] \textbf{Adding declared variables into a model.}
Before one can start to specify the
 objective function and constraints in a model,
 the variables, say {\tt X} and {\tt Y}, that  we have declared
 must be added to the {\tt ccp\_model} class  \texttt{mymodel} that
 we have created before.  This step is simply done via the command:
\begin{center}
	{\tt mymodel}.\textbf{add\_variable}({\tt X},{\tt Y});
\end{center}
Here {\tt mymodel} is a class object and \textbf{add\_variable} is a method in the class.

\end{description}

%
%

\subsection{Declaring the objective function}
\label{subsec-obj}

After creating the model {\tt mymodel}, declaring variables (say {\tt X} and {\tt Y})
and adding them into
{\tt mymodel}, we can  proceed to specify the objective function. Declaring an objective function requires the use of the functions (methods) \textbf{minimize} or \textbf{maximize}. There must be one and only one objective function
 in a model specification. In general, the objective function is specified
 through the sum or difference of the \textbf{inprod} function (inner product of two vectors or two matrices) which must have two input arguments in the form:
\textbf{inprod}({\tt C},{\tt X})
where {\tt X} must be a declared variable, and {\tt C} must be a constant vector or  matrix which is already available in the workspace and having the same dimension as {\tt X}. The input {\tt C} can also be a constant vector or matrix generated by some {\sc Matlab} built-in functions such as \texttt{speye(n,n)}.

Although we encourage users to specify an optimization problem in the standard form given in \eqref{eq-gensdp}, as a user-friendly interface, we also provide some extra functions to help users to specify the
objective function in a more natural way. We summarize these functions and their usages in
Table \ref{table: objective-func}.

\begin{table}[!h]
	\begin{tabular}{|c|p{11cm}|}
	\hline
	Function & Description \\
	\hline
	\textbf{inprod}({\tt C}, {\tt X}) & The inner product of a constant vector or matrix {\tt C} and variable {\tt X} of the same dimension.\\
	\hline
	\textbf{trace}({\tt X}) & The trace of a square matrix variable {\tt X}.\\
	\hline
	\textbf{sum}({\tt X}) & The sum of all elements of a vector or matrix variable {\tt X}.\\
	\hline
	\textbf{l1\_norm}({\tt X}) & The $\ell_1$ norm of a variable {\tt X}.\\
	\hline
	\textbf{l1\_norm}($\cA\ast${\tt X} +{\tt  b}) & The $\ell_1$ norm of an affine expression. {For the exact meaning of the expression ``{$\cA\ast${\tt X}}'', the reader can refer to \eqref{eq-AX}.}\\
	\hline
	\end{tabular}
	\caption{Supported functions for specifying the objective function in a model.}
	\label{table: objective-func}
\end{table}

For the class {\tt mymodel} created in Listing \ref{code: list1}, we can see that
the objective function
of  \eqref{EX1-interface} is specified via the command:
\begin{center}
\begin{footnotesize}
\begin{tt}
mymodel.minimize(trace(X1) + trace(X2) + sum(X3));
\end{tt}
\end{footnotesize}
\end{center}

\subsection{Adding affine constraints into the model}
\label{subsec-constraints}

Affine constraints can be specified and added into \texttt{mymodel} after the relevant variables have been declared.
This is done via the function (method) \texttt{add\_affine\_constraint}.
The following constraint types are supported in the interface:
\begin{itemize}
	\item
	Equality constraints \texttt{==}
	
	\item
	Less-or-equal inequality constraints \texttt{<=}
	\item
	Greater-or-equal  inequality constraints \texttt{>=}
\end{itemize}
where the expressions on both the left and right-hand sides of the operands must be affine expressions.
Strict inequalities \texttt{<} and \texttt{>} are \textbf{not} accepted.
Inequality and equality constraints are applied in an elementwise fashion, matching the behavior of {\sc Matlab} itself. For instance, if \texttt{U} and \texttt{X} are $m \times n$ matrices, then \texttt{X <= U} is interpreted as $mn$ (scalar) inequalities \texttt{X(i,j) <= U(i,j)} for all $i=1,\dots,m$, $j=1,\dots,n$. When one side is a scalar and the other side is a variable, that value is replicated; for instance, \texttt{X >= 0} is interpreted as \texttt{X(i,j) >= 0} for all $i=1,\dots,m$, $j=1,\dots,n$.

In general, affine constraints have the following form
\begin{equation}\label{Eq: Affine_General}
\cA_1 * {\tt X_1} + \cA_2*{\tt X_2} + \cdots + \cA_k*{\tt X_k}\; \mbox{\tt <=}\;
 (\mbox{\tt >=} \;\; {\rm or} \;\; \mbox{\tt ==}) \; b,
\end{equation}
where ${\tt X_1}, {\tt X_2}, \dots, {\tt X_k}$ are declared variables,  $b$ is a constant matrix or vector,
and $\cA_1, \cA_2, \dots, \cA_k$ are linear maps
whose descriptions will be given shortly.

Next, we illustrate how to add affine constraints into the model object \texttt{mymodel} in detail.

\subsubsection{General affine constraints}
In this section, we show users how to initialize the linear maps $\cA_1$, $\cA_2$, $\dots$, $\cA_k$ in \eqref{Eq: Affine_General}.
\begin{itemize}
	\item If $\cA_i = a_i$, is a scalar, then $a_i *{\tt X_i}$ has the same dimension as the {variable} ${\tt X_i}$.
	\item If ${\tt X_i}$ is an $n$-dimensional vector, then $\cA_i$ must be a $p \times n$ constant matrix, and $\cA_i *{\tt X_i}$ is in $\Re^{p}$.
	\item If ${\tt X_i}$ is an $m \times n$ ($n > 1$) matrix, then $\cA_i * {\tt X_i}$ is interpreted as
	a linear map such that
	\begin{eqnarray}
	\cA_i * {\tt X_i} = \left[ \begin{array}{c} \inprod{A^{(i)}_1}{{\tt X_i}}
	\\ \vdots \\ \inprod{A_p^{(i)}}{{\tt X_i}}
	\end{array}
	\right] \in \Re^p,
	\label{eq-AX}
	\end{eqnarray}
	where $A_1^{(i)},\ldots,A_p^{(i)}$ are given $m\times n$ constant matrices.
	In this case, $\cA_i$ is a $p \times 1$ constant cell array such that
	\begin{equation*}
	\mathcal{A}_i\{ j \} = A_j^{(i)}, \quad j=1,\dots,p.
	\end{equation*}
	
\end{itemize}

\subsubsection{Coordinate-wise affine constraints}

Although users can model coordinate-wise affine constraints in the general  form
given in \eqref{Eq: Affine_General}, we allow users to declare them in a more direct way as follows:
\begin{equation}\label{Eq: Affine_coordinate}
a_1*{\tt X_1}(i_1, j_1) + a_2*{\tt X_2}(i_2, j_2) + \dots + a_k*{\tt X_k}(i_k, j_k) \;
\mbox{\tt <=} \;  (\mbox{\tt >=}\;\;  {\rm or} \;\;\mbox{\tt ==}) \; b,
\end{equation}
where $a_1, a_2, \dots, a_k, b$ are scalars and ${\tt X_1}, {\tt X_2}, \dots, {\tt X_k}$ are declared variables.
The index pairs
$(i_1, j_1)$, $(i_2, j_2)$, $\dots$, $(i_k, j_k)$ extract the corresponding elements in the variables.
From Listing \ref{code: list1}, we can see how a constraint of the form \eqref{Eq: Affine_coordinate} is
added, i.e.,
\begin{center}	
{\tt mymodel}.add\_affine\_constraint($2*{\tt X_1}(2,3) + {\tt X_2}(4,2)  - {\tt X_3}(4) \;\mbox{\tt ==}\; 3$)
\end{center}

Our interface also allows users to handle multiple index pairs. For example, if we have a declared variable ${\tt X} \in \Re^{m \times n}$ and two index arrays
$$I = [i_1, i_2, \dots, i_k], \ \ J = [j_1, j_2, \dots, j_k],$$
where
$\max\{i_1, i_2, \dots, i_k\} \leq m $ and $ \max\{j_1, j_2, \dots, j_k\} \leq n$,
then ${\tt X}(I, J)$ is interpreted as
$${\tt X}(I, J) = \left[\begin{matrix}
{\tt X}(i_1, j_1)\\
{\tt X}(i_2, j_2)\\
\vdots\\
{\tt X}(i_k, j_k)\\
\end{matrix}\right] \in \Re^k.
$$
An example of such a usage can  be found in Listing \ref{code: ex_edm}.

\subsubsection{Element-wise multiplication}
In our interface, we also support element-wise multiplication $(.*)$ between a declared
variable {\tt X} and a constant matrix $A$ with the same dimension. Suppose
\begin{equation*}
{\tt X} = \left[
\begin{matrix}
{\tt X_{11}} & \cdots &{\tt X_{1n}}\\
\vdots  &\ddots &\vdots\\
{\tt X_{m1}}  & \cdots & {\tt X_{mn}}
\end{matrix}
\right], \quad  A = \left[\begin{matrix}
A_{11} & \cdots & A_{1n}\\
\vdots & \ddots & \vdots\\
A_{m1} & \cdots & A_{mn}
\end{matrix}\right].
\end{equation*}
Then $A.*{\tt X}$ is interpreted as
\begin{equation*}
A.*{\tt X} = \left[\begin{matrix}
A_{11}*{\tt X_{11}} & \cdots & A_{1n}*{\tt X_{1n}}\\
\vdots & \ddots & \vdots\\
A_{11}*{\tt X_{m1}} & \cdots & A_{mn}*{\tt X_{mn}}
\end{matrix}\right].
\end{equation*}

\subsubsection{Specifying affine constraints using predefined maps}
\label{subsec-constr_map}

For convenience, we also provide some predefined maps to help users to specify constraints in a more direct way. We summarize these maps and their usages in Table \ref{table: maps_affine_constr}.

\begin{table}[!h]
	\begin{tabular}{|c|p{8cm}|c|}
		\hline
		Function & Description & Dimension\\
		\hline
		\textbf{inprod}({\tt C}, {\tt X}) & The inner product of a constant vector or matrix {\tt C} and a variable {\tt X} of the same dimension. & $1\times 1$\\
		\hline
		\textbf{trace}({\tt X}) & The trace of a square matrix variable {\tt X}. & $1\times 1$\\
		\hline
		\textbf{sum}({\tt X}) & The sum of all elements of a vector or matrix variable {\tt X}. & $1\times 1$\\
		\hline
		\textbf{l1\_norm}({\tt X}) & The $\ell_1$ norm of a variable {\tt X}. & $1\times 1$\\
		\hline
		\textbf{l1\_norm}($\cA$*{\tt X} + b) & The $\ell_1$ norm of an affine expression.&$1\times 1$\\
		\hline
		\textbf{map\_diag}({\tt X}) & Extract the main diagonal of an $n\times n $ matrix
		variable {\tt X}. & $n \times 1$\\
		\hline
		\textbf{map\_svec({\tt X})} & For an $n\times n$ symmetric variable {\tt X}, it returns the corresponding symmetric vectorization of {\tt X}, as defined in \eqref{fun-svec}. & $\frac{n(n+1)}{2}\times 1$\\
		\hline
		\textbf{map\_vec}({\tt X}) & For a $m\times n$ matrix variable {\tt X}, it returns the vectorization of {\tt X}. & $mn \times 1$\\
		\hline
	\end{tabular}
	\caption{Supported predefined maps. }
	\label{table: maps_affine_constr}
\end{table}
\subsubsection{Chained constraints}

In our interface,
one can add chained inequalities into the created {\tt ccp\_model} \texttt{mymodel}. In general, chained affine constraints have the form
\begin{equation*}
\mbox{\tt L <=} \; \cA_1*{\tt X_1} + \cA_2*{\tt X_2} + \cdots + \cA_k*{\tt X_k}\; \mbox{\tt <= U},
\end{equation*}
where {\tt L} and {\tt U} are scalars or constant matrices with having the same dimensions as the affine expression in the middle.
As an example, one can add bound constraints for a declared variable {\tt X} via the command:
\begin{center}
\texttt{mymodel}.\texttt{add\_affine\_constraint}(\texttt{L <= X <= U});
\end{center}
It is important to note that in chained inequality constraints, the affine expression in the middle should only contain declared variables but not constants.

\subsection{Adding positive semidefinite constraints into the model}
\label{subsec-constr_psd}

Positive semidefinite constraints can be added into a previously created object
\texttt{mymodel} using the function (method) \texttt{add\_psd\_constraint}.  Such a constraint is
valid only for a declared symmetric variable or positive semidefinite variable.
In general, a positive semidefinite constraint has the form
\begin{equation}\label{constr: psd}
a_1*{\tt X_1} + a_2*{\tt X_2} + \cdots + a_k*{\tt X_k}\; \succeq \;{\tt G},
\end{equation}
where $a_1, a_2, \dots, a_k$ are scalars, and ${\tt X_1}$, ${\tt X_2}$, $\dots$, ${\tt X_k}$ are declared variables in symmetric matrix spaces or PSD cones, and {\tt G} is a constant symmetric matrix. Note that one can also have the version
 ``$\preceq$'' in \eqref{constr: psd}.
We can add \eqref{constr: psd} into \texttt{mymodel} as follows:
\begin{center}
	\texttt{mymodel.add\_psd\_constraint}($a_1*{\tt X_1} + \cdots + a_k*{\tt X_k}$ {\tt >=}  {\tt G})
\end{center}
Specially,
\begin{itemize}
	\item For a variable ${\tt X} \in \S^{n}$, one can use \texttt{mymodel.add\_psd\_constraint({\tt X}>=0)} to specify the constraint ${\tt X} \succeq 0$ or ${\tt X} \in \S^{n}_+$.
	\item For a variable ${\tt X} \in \S^{n}$ and a constant matrix ${\tt G} \in \S^{n}$. One can use \texttt{mymodel.add\_psd\_constraint({\tt X} >= G)} and \texttt{mymodel.add\_psd\_constraint({\tt X} <= G)} to specify the constraint ${\tt X} \succeq {\tt G}$ and ${\tt X} \preceq {\tt G}$, respectively.
\end{itemize}
Similar to affine constraints, one can also use chained positive semidefinite constraints together. For example,
for a variable ${\tt X} \in \S^n$ and two constant matrices ${\tt G1}, {\tt G2} \in \S^n$ $({\tt G1} \preceq {\tt G2})$, one can specify ${\tt G1} \preceq {\tt X} \preceq {\tt G2}$ as
\begin{center}
	\texttt{mymodel.add\_psd\_constraint(G1 <= X <= G2)};
\end{center}
\subsection{Setting parameters for \SDPNALP}
\label{subsec-setpara}

As described in Section \ref{subsec-para}, there are mainly nine parameters in the parameter structure array \texttt{OPTIONS}. To allow users to set these parameters freely, we provide the function (method) \texttt{setparameter} for such a
purpose. Parameters which are not specified are set to be the default values described in Section \ref{subsec-para}. Now, we describe the usage of \texttt{setparameter} in details.

{Assume} that we have created a \texttt{ccp\_model} class called
\texttt{mymodel}.  {Since} \texttt{setparameter} is a method in the \texttt{ccp\_model} class, so the usage of \texttt{setparameter} is simply
\begin{center}
	\texttt{{mymodel}.setparameter(`para\_name',value)}
\end{center}
In Table \ref{table: usage_setpara}, we summarize the parameters which can be set  in \texttt{setparameter}.
Note that users can set more than one parameters at a time. For example, one can use
\begin{center}
	\texttt{{mymodel}.setparameter(`\texttt{tol}', 1e-4, `maxiter', 2000)};
\end{center}
to set the parameters \texttt{tol = 1e-4} and \texttt{maxiter = 2000}.

\begin{table}[!h]
	\begin{center}
	\begin{tabular}{|c|l|c|}
		\hline
		Parameter Name & Usage & Default Value\\
		\hline
		\texttt{tol} & \texttt{{mymodel}.setparameter(`\texttt{tol}', value)} & 1e-6\\ \hline
	\texttt{maxiter} & \texttt{{mymodel}.setparameter(`\texttt{maxiter}', value)} & 20000\\ \hline
	\texttt{maxtime} & \texttt{{mymodel}.setparameter(`\texttt{maxtime}', value)} & 10000\\
		\hline
		\texttt{tolADM} & \texttt{{mymodel}.setparameter(`\texttt{tolADM}', value)} & 1e-4\\
		\hline
		\texttt{maxiterADM} & \texttt{{mymodel}.setparameter(`\texttt{maxiterADM}', value)} & 200\\
		\hline
		\texttt{printlevel} & \texttt{{mymodel}.setparameter(`\texttt{printlevel}', value)} & 1\\
		\hline
		\texttt{stopoption} & \texttt{{mymodel}.setparameter(`\texttt{stopoption}', value)} & 1\\
		\hline
		\texttt{AATsolve.method} & \texttt{{mymodel}.setparameter(`\texttt{AATsolve.method}', value)} & `direct'\\
		\hline
		\texttt{BBTsolve.method} & \texttt{{mymodel}.setparameter(`\texttt{BBTsolve.method}', value)} & `iterative'\\
		\hline
	\end{tabular}
    \end{center}
	\caption{Usage of \texttt{setparameter}.}
	\label{table: usage_setpara}
\end{table}

\subsection{Solving a model and extracting solutions}

After creating and initializing the class \texttt{mymodel}, one can call the method \textbf{solve}  to
solve the model as follow:
\begin{center}
	\texttt{{mymodel}.solve}
\end{center}
After solving the SDP problem, one can extract the optimal solutions using the function \textbf{get\_value}.
For example, if \texttt{{\tt X1}} is a declared variable, then
one can extract the optimal value of \texttt{{\tt X1}} by setting
\begin{center}
                          \texttt{get\_value({\tt X1})}
\end{center}
Note that the input of the function \texttt{get\_value} should be a declared variable.

%
%
%

\subsection{Further remarks on the interface}

Here we give some remarks to help users to input an SDP
 problem into our interface more efficiently.
\begin{itemize}
	\item If a variable must satisfy a conic constraint, it would be more efficient to specify the
	conic constraint when declaring the variable rather than declaring the variable and imposing the
	constraint separately. For example, it is better to use
	 \texttt{{X} = var\_nn(m,n)} to indicate that the variable ${\tt X}\in \Re^{m\times n}$ must be in the cone
	 $\Re^{m\times n}_+$ rather than separately declaring
	 {\tt {X} = var\_free(m,n)} followed by setting
	 \begin{center}
	 {\tt {mymodel}.add\_affine\_constraint({\tt X} >= 0);}
	\end{center}
	Similarly, if a square matrix variable ${\tt Y}\in \S^n$ must
	satisfy the conic constraint that ${\tt Y}\in\S^n_+$, then it
	is better to declare it as \texttt{{Y} = var\_sdp(n,n)} rather than
	 separately declaring \texttt{{Y} = var\_free(n,n)} followed by setting
	 \begin{center}
	 {\tt {mymodel}.add\_psd\_constraint({\tt Y} >= 0);}
	\end{center}
	The latter option is not preferred because we have to introduce extra constraints.
	\item When there is a large number of affine constraints, specifying them using a loop
	in {\sc Matlab} is generally time consuming. To make the task more efficient, if possible, always try to model the problem using our predefined functions
\end{itemize}

\section{Examples on building SDP models using our interface}
\label{sec-interface-examples}
To solve SDP problems using \SDPNALP\!\!, the user must input the problem data corresponding to the  form in (P). The file {\tt SDPNALplusDemo.m} contains a few examples to illustrate how to
 generate the data of an SDP problem in the required format. Here we will present a few of those
 examples in detail. Note that the user can also store the problem data in either the SDPA or SeDuMi format, and then use the m-files to read {\tt sdpa.m} or {\tt sedumi.m} to convert the data for \SDPNALP\!\!.

 We also illustrate how the SDP problems can be coded
 using our basic interface.

\subsection{SDPs arsing from the nearest correlation matrix problems}
\label{sec-NCM-interface}

To obtain a valid nearest correlation matrix (NCM) from a given incomplete sample
correlation matrix $G\in\S^n$, one version of the NCM problem is to consider solving
the following SDP:
\begin{eqnarray*}
\mbox{(NCM)}  && \min \Big\{ \norm{H\circ(X-G)}_1
  \;\mid\; {\rm diag} (X) = e, \quad X \in \S^n_+ \Big\},
\end{eqnarray*}
where $H\in \S^n$ is a nonnegative weight matrix and ``$\circ$" denotes the elementwise
product. Here for any $M\in \S^n$, $\norm{M}_1 = \sum_{i,j=1}^n |M_{ij}|$.

In order to express (NCM) in the form given in (P), we first write
$$
{\rm svec}(X) - {\rm svec}(G) = x_+ -  x_-,
$$
where $ x_+$ and $x_-$ are two nonnegative vectors in $\Re^{\bar{n}}$ ($\bar{n} = n(n + 1)/2$). Then  (NCM) can be reformulated as the following SDP with $m = n+\bar{n}$ equality constraints:
\begin{eqnarray}
\begin{array}{cl}
\min & \inprod{{\rm svec}(H)}{x_+} +  \inprod{{\rm svec}(H)}{x_-} \\[4pt]
{\rm s.t.} & {\rm diag}(X) \qquad\qquad\;\;\; = e, \\[4pt]
& {\rm svec}(X) - x_+  +  x_- = {\rm svec}(G), \quad
X \in \S^n_+, \; x_+,  x_- \in \Re^{\bar{n}}_+.
\end{array}
\label{eq-NCM-2}
\end{eqnarray}
Given $G,H\in\S^n$, the SDP data for the above problem can be coded for \SDPNALP as follows.

\begin{footnotesize}
\begin{tt}
\begin{lstlisting}[caption = {Generating the \SDPNALP data for the NCM problem \eqref{eq-NCM-2}.},frame = single]
blk{1,1} = 's'; blk{1,2} = n;
n2 = n*(n+1)/2;
II = speye(n2); hh = svec(blk(1,:),H);

for k=1:n; Acell{k} = spconvert([k,k,1;n,n,0]); end
Atmp = svec(blk(1,:),Acell,1);
At{1,1} = [Atmp{1}, II];
At{2,1} = [sparse(n,n2), sparse(n,n2); -II, II]';

b = [ones(n,1); svec(blk(1,:),G)];
C{1,1} = sparse(n,n); C{2,1} = [hh; hh];
\end{lstlisting}
\end{tt}
\end{footnotesize}
For more details, see the m-file {\tt NCM.m} in the subdirectory {\tt /util}.

Next, we show how to use our interface to solve the nearest correlation matrix problem
(NCM).
Given a data matrix  $G\in \S^n$, we can
solve the corresponding NCM problem using our interface as follows.

\begin{footnotesize}
\begin{tt}
\begin{lstlisting}[caption = {Solving a NCM problem with our interface.}
\label{code:ex_ncm},frame = single]
 n = 100;
 G = randn(n,n);  G = 0.5*(G + G');
 H = rand(n);  H = 0.5*(H+H');
 model = ccp_model('Example_NCM');
    X = var_sdp(n,n);
    model.add_variable(X);
    model.minimize(l1_norm(H.*X - H.*G));
    model.add_affine_constraint(map_diag(X) == ones(n,1));
    model.setparameter('tol', 1e-6, 'maxiter', 2000);
 model.solve;
 Xval = get_value(X);
 dualinfo = get_dualinfo(model);
\end{lstlisting}
\end{tt}
\end{footnotesize}

The last two lines in Listing \ref{code:ex_ncm} illustrate how we can
extract the numerical value of the variable $X$ and also the corresponding
dual variables.
Observe that with the help of our interface, users can
input the problem into our solver very easily; see \texttt{Example\_NCM.m} for more details.
\subsection{SDP relaxations of the maximum stable set problems}

Let $G$ be an undirected graph with $n$ nodes and edge set $\mathcal{E}$.
Its stability number, $\alpha(G)$, is the cardinality of a maximal stable set of $G$, and
it can be expressed as
\[
\alpha(G) := \max \{ e^Tx \,:\, x_ix_j =0, (i,j)\in \mathcal{E}, x \in \{0,1\}^n\},
\]
where $e \in \Re^n$ is the vector of all ones. It is known that computing $\alpha(G)$ is
NP-hard. But an upper bound $\theta(G)$, known as the
Lov\'{a}sz theta number \cite{LovaszTheta}, can be computed as the optimal value of the following SDP problem:
\begin{eqnarray}\label{eq-theta}
\theta(G) &:= & \max  \Big\{ \inprod{ee^T}{X}
\;\Big|\; \inprod{E^{ij}}{X} =0\;  \forall \; (i,j) \in \mathcal{E}, \;
   \inprod{I}{X} =1, \; X \in \S^n_+ \Big\},
\end{eqnarray}
where $E^{ij}=e_ie_j^{\T} + e_je_i^{\T}$ and $e_i$ denotes the $i$th standard unit vector
of $\Re^n$. One can further tighten the upper bound to get
$\alpha(G) \leq \theta_+(G) \leq \theta(G)$, where
\begin{eqnarray}
\label{eq-theta-plus}
\theta_+(G) &:= & \max  \Big\{ \inprod{ee^T}{X}
\;\Big|\; \inprod{E^{ij}}{X} =0\;  \forall \; (i,j) \in \mathcal{E}, \;
   \inprod{I}{X} =1, \; X \in \S^n_+,\; X \geq 0 \Big\}.
\end{eqnarray}
In the subdirectory {\tt /datafiles} of \SDPNALP\!\!,  we provide a few SDP problems
with data stored in the
in SDPA or SeDuMi format,
arising from computing $\theta(G)$ for a few graph instances.
The segment below illustrates how one can solve the
SDP problem, {\tt theta8.dat-s}, to compute $\theta_+(G)$:
\begin{small}
\begin{verbatim}
>> [blk,At,C,b] = read_sdpa('theta8.dat-s');
>> L = 0;
>> [obj,X,s,y,S,Z,ybar,v,info,runhist] = sdpnalplus(blk,At,C,b,L);
\end{verbatim}
\end{small}
To compute $\theta(G)$, one can simply set {\tt L = []} to indicate that
there is no lower bound constraint  on $X$.
In Listing \ref{code:ex-theta-plus}, we illustrate how to use
our interface to solve the $\theta_+$ problem \eqref{eq-theta-plus}.

\begin{tt}
\begin{footnotesize}
\begin{lstlisting}[caption = {Solving the $\theta_+$ problem
\eqref{eq-theta-plus} using our interface.}\label{code:ex-theta-plus},frame = single,
backgroundcolor=\color{grey}]
 load theta6.mat
 [IE,JE] = find(triu(G,1));
 n = length(G);
 model = ccp_model('Example_theta');
     X = var_sdp(n,n);
     model.add_variable(X);
     model.maximize(sum(X));
     model.add_affine_constraint(trace(X) == 1);
     model.add_affine_constraint(X(IE,JE) == 0);
     model.add_affine_constraint(X >= 0);
 model.solve;
\end{lstlisting}
\end{footnotesize}
\end{tt}

\subsection{SDPs arising from the frequency assignment problems}
\label{sec-Interfaceex-FAPs}

Given a network represented by a graph $G$ with $n$ nodes and an edge set $\cE$ together with an edge-weight matrix $W$, a certain type of frequency assignment problem on $G$  can be relaxed into the following SDP (see \cite[eq. (5)]{burer2003computational}):
\begin{eqnarray} \label{eq-fap}
\begin{array}{rl}
\mbox{(FAP)} \quad
\max &  \inprod{(\frac{k-1}{2k})\mathcal{L}(G,W)-\frac{1}{2}{\rm Diag}(We)}{ X} \\[5pt]
{\rm s.t.} & {\rm diag}(X) = e, \quad X \in \S^n_+,
\\[3pt]
& \inprod{-E^{ij}}{ X} = 2/(k-1) \quad \forall\; (i,j)  \in \cU \subseteq \cE,
\\[3pt]
& \inprod{-E^{ij}}{ X} \leq 2/(k-1) \quad \forall\;  (i,j) \in \cE \setminus \cU,
\end{array}
\end{eqnarray}
where $k>1$ is a given integer, $\cU$ is a given subset of $\cE$, $\mathcal{L}((G,W) := {\rm Diag}(We)-W$ is the Laplacian matrix, $E^{ij}=e_ie_j^T + e_je_i^T$.
Note that (\ref{eq-fap}) is equivalent to
\begin{eqnarray} \label{eq-fap-sdpp}
\begin{array}{ll}
\max & \inprod{(\frac{k-1}{2k})\mathcal{L}(G,W)-\frac{1}{2}{\rm Diag}(We)}{ X}
\\[5pt]
{\rm s.t.} & {\rm diag}(X) = e, \;\; X  \in \S^n_+,\;\; L \leq X \leq U,
\end{array}
\end{eqnarray}
where
\begin{eqnarray*}
	L_{ij} =
	\begin{cases}
		-\frac{1}{k-1} &  \forall (i,j) \in \cE, \\
		- \infty & {\rm otherwise},
	\end{cases}
	\quad
	U_{ij} =
	\begin{cases}
		-\frac{1}{k-1} &  \forall (i,j) \in \cU, \\
		 \infty & {\rm otherwise}.
	\end{cases}
\end{eqnarray*}

Next, we show how to use our interface to solve the SDP problem \eqref{eq-fap}.
Assume that we have already computed the constant matrix $C := (\frac{k-1}{2k})\mathcal{L}(G,W)-\frac{1}{2}{\rm Diag}(We)$ and saved it as \texttt{C} in the current workspace.
Suppose
\texttt{IU}, \texttt{JU} are two single column arrays storing the index pairs $(i,j)$ corresponding to $\cU$, and
\texttt{IE}, \texttt{JE} are two single column arrays storing the index pairs $(i,j)$ corresponding to $\cE$.
Assume that {\tt IU}, {\tt JU}, {\tt IE}, {\tt JE}, {\tt n}, {\tt kpara} are already stored in
the current workspace.
We can build the \texttt{ccp\_model} for \eqref{eq-fap} using our interface as follows.
More details can be seen in \texttt{Example\_FAP.m}.
\begin{tt}
\begin{footnotesize}
\begin{lstlisting}[caption = {Solving the FAP \eqref{eq-fap} using our interface.}\label{code: ex_FAP},frame = single]
model = ccp_model('Example_FAP');
   X = var_sdp(n,n);
   model.add_variable(X);
   model.maximize(inprod(C,X));
   model.add_affine_constraint(map_diag(X) == ones(n,1));
   const = -1/(kpara-1);
   model.add_affine_constraint(X(IU,JU) == const);
   model.add_affine_constraint(X(IE,JE) >= const);
model.solve;
\end{lstlisting}
\end{footnotesize}
\end{tt}

One can also solve (FAP) using the equivalent formulation
specified in \eqref{eq-fap-sdpp}.
Assume that the matrices {\tt L}, {\tt U}, {\tt C} and {\tt n} have been computed in the current workspace,
we can input the SDP problem \eqref{eq-fap-sdpp}  into our interface based on the above equivalent form
as follows.
\begin{footnotesize}
\begin{tt}
\begin{lstlisting}[caption = {Solving the reformulated FAP \eqref{eq-fap-sdpp}.}\label{code: ex_fap_refor},frame = single]
 model = sdp_model('Example_FAP2');
   X = var_sdp(n,n);
   model.add_variable(X);
   model.maximize(inprod(C,X));
   model.add_affine_constraint(map_diag(X) == ones(n,1));
   model.add_affine_constraint(L <= X <= U);
 model.solve;
\end{lstlisting}
\end{tt}
\end{footnotesize}

\subsection{SDPs arising from Euclidean distance matrix problems}

Consider a given undirected graph $G$ with $n$ nodes and edge set $\cE$.
Let $D=(d_{ij})\in \S^n$ be a matrix whose elements are such that
$d_{ij} > 0$ if $(i, j) \in \cE$, and $d_{ij} = 0$ if $(i, j) \not\in \cE$.
We seek points $x_1, x_2, \ldots, x_n$ in $\Re^d$ such that $\norm{x_i-x_j}$ is
as close as possible to $d_{ij}$ for all $(i, j) \in \cE$. In particular, one
may consider minimizing the $L_1$-error as follows:
\begin{eqnarray*}
& \min \Big\{ \sum_{(i,j)\in \cE}  |d^2_{ij} - \|x_i -x_j\|^2|  -\frac{\alpha}{2n}
\sum_{i,j=1}^n\norm{x_i-x_j}^2  \,\mid\,
 \sum_{i=1}^n x_i= 0, \; x_1,\ldots,x_n \in \Re^d
\Big\},&
\end{eqnarray*}
where the equality constraint is introduced to put the center of mass of the points at the origin.
The second term in the objective function is introduced to achieve the effect of spreading out the points instead of crowding together, and $\alpha$ is a given nonnegative parameter.
Let $X=[x_1, \ldots, x_n] \in \Re^{d \times n}$. Then $\|x_i- x_j\|^2 = e^T_{ij} X^T X e_{ij}$, where $e_{ij} = e_i-e_j$. The above nonconvex problem can be rewritten as (for more details, see \cite{LeungToh2009}):
\begin{eqnarray*}\label{EDM2}
\begin{array}{rl}
\min & \Big\{ \sum_{(i,j) \in \cE} |d^2_{ij} - \inprod{e_{ij}e_{ij}^T}{Y} | -\alpha\inprod{I}{Y}\,\mid\,
  \inprod{E}{Y} = 0, \;
 Y = X^T X, \; X \in \Re^{d \times n}
\Big\}.
 \end{array}
\end{eqnarray*}
By relaxing the nonconvex constraint $Y=X^TX$ to $Y\in \S^n_+$, we obtain the
following SDP problem:
\begin{eqnarray}\label{EDM-SDP}
\begin{array}{rl}
\min & \sum_{(i,j) \in \cE} \; x^+_{ij}+x^-_{ij}  -\alpha\inprod{I}{Y} \\[5pt]
 {\rm s.t.} &  \inprod{e_{ij}e_{ij}^T}{Y} - x^+_{ij} + x^-_{ij} =  d^2_{ij}  \; \forall \; (i,j) \in \cE,
  \\[6pt]
  &\inprod{E}{Y} = 0,
 \\[5pt]
 &  Y \in \S^n_+, \; x^+_{ij}, x^-_{ij} \geq 0\; \forall \; (i,j) \in \cE.
\end{array}
\end{eqnarray}
Note that the number of the equality constraints in (\ref{EDM-SDP}) is $ |\cE|+1$,
and that the problem does not satisfy the Slater's condition because of the constraint
$\inprod{E}{Y} = 0$.
The problem (\ref{EDM-SDP}) is typically highly degenerate
and the optimal solution is not unique, which may result
 in high sensitivity to small perturbations in
the data matrix $D$. Hence, the problem (\ref{EDM-SDP}) can usually
only be solved by \SDPNALP to a
 moderate accuracy tolerance, say ${\tt OPTIONS.tol}= 10^{-4}$.
Given the data matrix $D\in\S^n$,
and let $m =|\cE|$,
the SDP data for \eqref{EDM-SDP}
can be  coded  as follows:

\begin{footnotesize}
\begin{tt}
\begin{lstlisting}[caption = {Generating the \SDPNALP data for the EDM problem \eqref{EDM-SDP}.},frame = single]
blk{1,1} = 's'; blk{1,2} = n;
Acell = cell(1,m+1); b = zeros(m+1,1); cnt = 0;
for i = 1:n
   for j = 1:n
      if (D(i,j) ~= 0)
         cnt = cnt + 1;
         Acell{cnt} = spconvert([i,i,1; i,j,-1; j,i,-1; j,j,1; n,n,0]);
         b(cnt) = D(i,j)^2;
      end
   end
end
Acell{m+1} = ones(n);
At(1) = svec(blk(1,:),Acell); C{1,1} = -alpha*speye(n,n);
blk{2,1} = 'l'; blk{2,2} = 2*m;
At{2,1} = [-speye(m), speye(m); sparse(1,2*m)]'; C{2,1} = ones(2*m,1);
\end{lstlisting}
\end{tt}
\end{footnotesize}

Next, we show how to solve the EDM problem \eqref{EDM-SDP}  using our interface.
Assume that we have generated the data matrix
$D\in \S^n$ such that $D_{ij} = d_{ij}$ for all $(i,j)\in \cE$, and stored it in \texttt{data\_randEDM.mat}
together with a given $\alpha$ . As mentioned above, we set the accuracy tolerance to solve the problem as \texttt{1e-4}. Now we can input the SDP problem into our interface as follows.
\begin{footnotesize}
\begin{tt}
\begin{lstlisting}[caption = {Solving the EDM problem \eqref{EDM-SDP} using our interface.}\label{code: ex_edm},frame = single]
 load data_randEDM;
 [ID, JD, val] = find(D);
 dd = val.^2;
 n1 = length(D);
 n2 = length(ID);

 model = ccp_model('Example_EDM');
   X1 = var_nn(n2,1);
   X2 = var_nn(n2,1);
   Y = var_sdp(n1,n1);
   model.add_variable(X1,X2,Y);
   model.minimize(sum(X1) + sum(X2) - alpha*trace(Y));
   model.add_affine_constraint(Y(ID,ID)+Y(JD,JD)-Y(ID,JD)-Y(JD,ID) -X1 +X2 == dd);
   model.add_affine_constraint(sum(Y) == 0);
   model.setparameter('tol', 1e-4, 'maxiter', 2000);
 model.solve;
\end{lstlisting}
\end{tt}
\end{footnotesize}

\subsection{SDPs arising from quadratic assignment problems}

Let $\Pi$ be the set of $n\times n$ permutation matrices. Given matrices $A, B \in \Re^{n\times n}$, the
associated quadratic assignment problem (QAP)  is given by
\begin{equation}\label{eq: QAP}
v^{*}_{QAP} : = \min\{\inprod{X}{AXB} : X \in \Pi\}.
\end{equation}
For a matrix $X = [x_1, \dots, x_n] \in \Re^{n\times n}$, we will identify it with the $n^2$-dimensional
vector $x = [x_1; \dots; x_n]$. For a matrix $Y \in \Re^{n^2\times n^2}$, we let $Y^{ij}$ be the $n\times n$ block corresponding to $x_ix_j^T$ in the matrix $xx^T$. It is shown in \cite{povh2009copositive} that $v^{*}_{QAP}$ is bounded below by the following number:
\begin{equation}\label{Eq: SDP_QAP}
\begin{array}{lll}
v := & \min & \inprod{B\otimes A}{Y}\\[5pt]
& s.t. & \sum_{i=1}^{n} Y^{ii} = I,\quad \inprod{I}{Y^{ij}} = \delta_{ij}  \; \; \forall \;\; 1\leq i\leq j \leq n,\\[5pt]
&& \inprod{E}{Y^{ij}} = 1, \; \; \forall\;\; 1\leq i\leq j\leq n,\\[5pt]
&& Y \succeq 0,\;\; Y \geq 0,
\end{array}
\end{equation}
where $E$ is the matrix of ones, and $\delta_{ij} = 1$ if $i=j$, and $0$ otherwise.
Note that there are $3n(n+1)/2$ equality constraints in \eqref{Eq: SDP_QAP}. But two of them are actually redundant, and we remove them when solving the standard SDP generated from \eqref{Eq: SDP_QAP}.

Now, we show an example of solving the SDP relaxation of the QAP problem \texttt{'chr12a'} via our interface.
\begin{footnotesize}
\begin{tt}
\begin{lstlisting}[caption = {Solving the SDP relaxation of a QAP with our interface.}\label{code: ex_qap},frame = single]
 problem_name = 'chr12a';
 [A, B] = qapread(strcat(problem_name, '.dat'));
 %% Construct C
 Ascale = max(1, norm(A, 'fro'));
 Bscale = max(1, norm(B, 'fro'));
 A = A/Ascale; B = B/Bscale;
 C = kron(B, A); C = 0.5*(C + C');
 nn = length(C);
 n = length(A);

 model = ccp_model(problem_name);
   Y = var_sdp(nn, nn);
   model.add_variable(Y);
   model.minimize(inprod(C, Y));
   model.add_affine_constraint(Y >= 0);
   II = speye(n); EE = ones(n);
   for i = 1:n-1
      for j = i:n
         Eij = sparse(i,j,1,n,n);
         if (i==j) const = 1; else, const = 0; end
         model.add_affine_constraint(inprod(kron(II,Eij), Y) == const);
         model.add_affine_constraint(inprod(kron(Eij,II), Y) == const);
         model.add_affine_constraint(inprod(kron(Eij,EE), Y) == 1);
      end
   end
   model.add_affine_constraint(inprod(kron(II,sparse(n,n,1,n,n)), Y) == 1);
   model.setparameter('maxiter', 5000);
 model.solve;
\end{lstlisting}
\end{tt}
\end{footnotesize}

\subsection{Comparison of our basic interface with CVX and YALMIP}

 As mentioned in the Introduction, our new interface
 is motivated by the need to
facilitate the modeling of an SDP problem for \SDPNALP
 to directly accept inequality constraints of the form
``$l \leq {\cal B}(X) \leq u$'', and bound constraints of the form
``$L\leq X\leq U$'' in addition to equality constraints
of the form ``$\cA (X) = b$''.

For the interfaces CVX \cite{CVX} and YALMIP \cite{YALMIP},  one will need to first reformulate
a problem with the above mentioned inequality constraints
  into the standard primal SDP form (for interior-point solvers) by converting
the inequality constraints into equality constraints through introducing
 extra nonnegative variables as follows:
\begin{eqnarray*}
& {\cal B}(X) - s^{(1)} = l, \; {\cal B}(X) + s^{(2)} = u,\;
  X- X^{(1)} = L,\; X+ X^{(2)} = U,
  \;&
  \\[5pt]
 & s^{(1)} \geq 0, \; s^{(2)} \geq 0, \; X^{(1)} \geq 0, \; X^{(2)} \geq 0.
 &
\end{eqnarray*}
The above conversion not only will add significant overheads when generating
the SDP data in CVX or  YALMIP, a much more serious computational issue
is that it has created a large number of additional equality
constraints in the formulation which would cause huge computational inefficiency
when solving the problem. Moreover, the large number of additional
equality constraints introduced  will likely make the SDP solver to
encounter various numerical difficulties when solving the resulting SDP problem.

In Table \ref{table-1}, we present the relevant information for the
SDP data generated by various interfaces for the QAP problem \eqref{Eq: SDP_QAP}  with
matrices $A,B$ of
dimensions $n\times n$. As one can observe, CVX took an exceeding long
time to generate the data compared to YALMIP and \SDPNALP\!\!.
When the problem dimension $n$ becomes larger, the ratio
of the times taken by YALMIP and \SDPNALP to generate the data
also grows larger, and the ratio is more than 13 for $n=20$.
More alarmingly, the number of equality constraints generated
by CVX or YALMIP is exceedingly large.
For $n=20$,
the ratio of the number of equality constraints generated by
YALMIP and \SDPNALP  is more than $255\, (\approx 160400/628)$ times.
Such a huge number of equality constraints generated by CVX or YALMIP
is fatal for the computational efficiency of interior-point solvers, and also disadvantageous for
\SDPNALP\!\!.

\begin{table}[h]
\begin{footnotesize}
\caption{Time taken (in seconds) to generate the SDP data
(and the corresponding problem sizes)
by various interfaces for the QAP problem \eqref{Eq: SDP_QAP} with
matrices $A,B$  of
dimension $n\times n$. Here
$m$ is the final number of equality constraints in the generated SDP data, ${\rm sblk}$ is the dimension
of the positive semidefinite matrix block, ${\rm lblk}$ is the dimension of the
nonnegative vector, ${\rm ublk}$ is the dimension of the
unrestricted vector.}
\label{table-1}
\begin{center}
\begin{tabular}{|c|c|c|c|} \hline
\mc{1}{|c|}{$n$} &\mc{1}{|c|}{CVX} &\mc{1}{|c|}{YALMIP} &\mc{1}{|c|}{\SDPNALP}
\\[5pt] \hline
10 & 9.22  & 2.55   & 0.49 \\
   & $\begin{array}{l} m=5213 \\ {\rm sblk}=100,\\ {\rm lblk}=5050\end{array}$
   &$\begin{array}{l} m=10100 \\ {\rm sblk}=100\\ {\rm lblk}=5050\\ {\rm ublk}=10163\end{array}$
   & $\begin{array}{l} m=163 \\ {\rm sblk}=100,\\ {\rm lblk}=5050\end{array}$
\\[5pt] \hline
15 &448  & 3.52 & 0.67  \\
& $\begin{array}{l} m=25783 \\ {\rm sblk}=225,\\ {\rm lblk}=25425\end{array}$
 &$\begin{array}{l} m=50850 \\ {\rm sblk}=225\\ {\rm lblk}=25425\\ {\rm ublk}=50983\end{array}$
 & $\begin{array}{l} m=358 \\ {\rm sblk}=225\\ {\rm lblk}=25425\end{array}$
 \\ \hline
 20 &  &9.86  &0.73  \\
 &\mbox{took too long to run}
 &$\begin{array}{l} m=160400 \\ {\rm sblk}=400\\ {\rm lblk}=80200\\ {\rm ublk}=160628\end{array}$
 & $\begin{array}{l} m=628 \\ {\rm sblk}=400\\ {\rm lblk}=80200\end{array}$
 \\ \hline
\end{tabular}
\end{center}
\end{footnotesize}
\end{table}

\section{Summary of the numerical performance of \SDPNALP}
\label{sec-num}

We have tested our solver \SDPNALP on
745 SDP instances arising from various sources, namely,
\begin{enumerate}
\item 65 instances of DNN (doubly nonnegative) relaxation of maximum
stable set problems from \cite{toh2004,Sloane,DIMACS};
\item 14 instances of SDP relaxation of frequency assignment problems (FAPs) \cite{FAP};
\item 94 instances of DNN relaxation of quadratic assignment problems (QAPs) \cite{QAPLIB};
\item 165 instances of DNN relaxation of  binary
quadratic integer programming (BIQ) problems \cite{BiqMac};
\item 120 instances of DNN relaxation of clustering problems \cite{Peng2007clustering};
\item 165 instances of DNN relaxation of BIQ problems with additional valid inequalities \cite{ADMM3c};

\item 65 instances of SDP relaxation of maximum stable set problems  from \cite{toh2004,Sloane,DIMACS};
\item 57 instances of SDP relaxation of best rank-one tensor
approximation  problems \cite{NWang2014}.
\end{enumerate}
In total there are 623 SDP problems with simple polyhedral bound constraints
on the matrix variable in addition to other linear constraints,
and 122 standard SDP problems.
The complete numerical results
are available at
\begin{center}
\url{http://www.math.nus.edu.sg/\~mattohkc/papers/SDPNALPtable-2017-Dec-18.pdf}
\end{center}
Note that the results are obtained on
a desktop computer having the following
specification: Intel Xeon CPU E5-2680v3 @2.50 GHz with 12 cores, and 128GB of RAM.
The extensive numerical experiments
 show that our \SDPNALP solver  is quite efficient and robust,
 in that it is able to solve 98.9\% of the 745  instances of  SDP problems
 arising from various applications listed above to the accuracy of less than
 $1.5\times 10^{-6}$ in the relative KKT residual $\eta$ defined in  \eqref{eq-eta}.

\begin{figure}[!t]
	\begin{center}
		\includegraphics[width=0.8\textwidth]{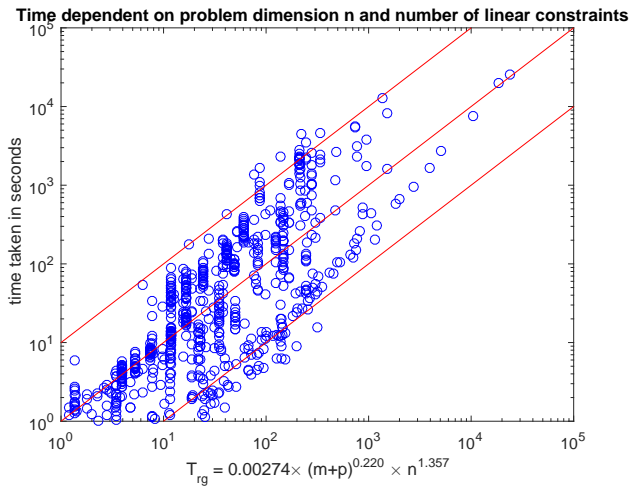}
			\vskip-5cm
		\caption{Time $T$ taken to solve 707 SDP instances versus the times estimated
		based on regression $T_{\rm rg} = 0.00274\, (m+p)^{0.220} \, n^{1.357}$. }
		\label{fig-1}
	\end{center}
\end{figure}

In Figure \ref{fig-1}, we plot the time $T$  taken to solve a
subset of 707 tested instances (with
computation time of over one second each)
versus the estimated times $T_{\rm rg} = 0.00274\, (m+p)^{0.220} \, n^{1.357}$,
obtained based on the regression $\log_{10} (T) \approx
\log_{10}(\kappa) + \alpha \log_{10} (m+p) + \beta \log_{10}(n).$
From the graph, one can observe that $T_{\rm rg}$ can estimate the actual
time taken to within a factor of about $20$ for a given $(m+p,n)$.
If we contrast the dependent of $T_{rg}$ on $(m+p,n)$ with the
$O((m+p)^2n^2) + O((m+p) n^3) + O((m+p)^3)$ time complexity
in an interior-point method such as those implemented in SDPT3 or SeDuMi,
then we can immediately observe that the time complexity of \SDPNALP
is much better. In particular, the {dependence on the number of}
linear constraints is only $(m+p)^{0.22}$ for a given matrix dimension $n$.
This also explains why our solver can be so efficient in solving an SDP problem
with a large number of linear constraints.

In Table \ref{table-2}, we give a summary of the numerical results obtained
for the subset of 707 SDP problems mentioned in the last paragraph.
Note that in the table, $m+p$ is the total number of linear constraints
as specified by $\cA$ and $\cB$. The simple polyhedral bound constraints
on the matrix variable are not counted in $m+p$. Thus even if $m+p$ is a
modest number, say less than 1000, the number of actual
polyhedral constraints in the problem can still be large.
Observe that across each row in the table, the average time taken
to solve the problems with different  number of linear constraints does
not depend strongly on $m+p$.
However, across each column in the table, the {dependence}
of the average time taken to solve the problems on the matrix
dimension $n$ is more significant, but it is
still much weaker than the cubic exponent dependent on the matrix dimension.

\begin{table}[!h]
\begin{footnotesize}
\caption{Summary of numerical results obtained by \SDPNALP in solving
707 SDP problems (each with the computation time of more than one second).
In each cell, the first number is the number of problems solved, and the
second number is the average time taken to solve the problems. Here $K$ means
a thousand.}
\label{table-2}
\begin{center}
\begin{tabular}{|c|c|c|c|c|c|c|c|} \hline
\mc{1}{|c|}{$m+p$} &\mc{1}{|c|}{$\leq 1K$} &\mc{1}{|c|}{$(1K,4K]$}
&\mc{1}{|c|}{$ (4K,16K]$}  &\mc{1}{|c|}{$ (16K,64K]$}
&\mc{1}{|c|}{$(64K,256K]$} &\mc{1}{|c|}{$(256K,1024K]$}
 &\mc{1}{|c|}{$> 1024K$}
\\[5pt] \hline

$n\leq 100$ &  $\begin{array}{l}  36 \\ 2.04\end{array}$ &  $\begin{array}{l}  16 \\ 2.16\end{array}$ &  $\begin{array}{l}  11 \\ 8.74\end{array}$ &  &  &  & 
\\[5pt] \hline
$100 < n\leq 200$ &  $\begin{array}{l} 101 \\ 6.29\end{array}$ &  $\begin{array}{l}   4 \\ 9.48\end{array}$ &  $\begin{array}{l}  43 \\ 25.67\end{array}$ &  $\begin{array}{l}  43 \\ 61.57\end{array}$ &  &  & 
\\[5pt] \hline
$200 < n\leq 400$ &  $\begin{array}{l} 127 \\ 29.17\end{array}$ &  $\begin{array}{l}   8 \\ 14.07\end{array}$ &  $\begin{array}{l}  15 \\ 3.88\end{array}$ &  $\begin{array}{l}  40 \\ 95.25\end{array}$ &  $\begin{array}{l}  20 \\ 269.95\end{array}$ &  & 
\\[5pt] \hline
$400 < n\leq 800$ &  $\begin{array}{l}  44 \\ 134.97\end{array}$ &  $\begin{array}{l}  14 \\ 603.19\end{array}$ &  $\begin{array}{l}   8 \\ 55.52\end{array}$ &  $\begin{array}{l}  16 \\ 34.75\end{array}$ &  $\begin{array}{l}  14 \\ 22.62\end{array}$ &  $\begin{array}{l}  15 \\ 2067.98\end{array}$ & 
\\[5pt] \hline
$800 < n\leq 1600$ &  $\begin{array}{l}  10 \\ 168.62\end{array}$ &  $\begin{array}{l}  56 \\ 871.80\end{array}$ &  $\begin{array}{l}  22 \\ 496.73\end{array}$ &  $\begin{array}{l}   7 \\ 115.58\end{array}$ &  $\begin{array}{l}  10 \\ 250.14\end{array}$ &  $\begin{array}{l}   5 \\ 172.92\end{array}$ & 
\\[5pt] \hline
$1600 < n\leq 3200$ &  &  $\begin{array}{l}   1 \\ 2672.64\end{array}$ &  &  $\begin{array}{l}   8 \\ 2867.77\end{array}$ &  $\begin{array}{l}   1 \\ 439.11\end{array}$ &  $\begin{array}{l}   5 \\ 2291.35\end{array}$ &  $\begin{array}{l}   1 \\ 966.32\end{array}$
\\[5pt] \hline
$n > 3200$ &  &  &  $\begin{array}{l}   1 \\ 12817.94\end{array}$ &  &  &  &  $\begin{array}{l}   5 \\ 11512.53\end{array}$
\\[5pt] \hline

\end{tabular}
\end{center}
\end{footnotesize}
\end{table}

\section{Conclusion and future works}
\label{sec-conclusion}

\SDPNALP  is   designed to be a general purpose software for solving
large scale SDP problems with bound constraints as well as having a
large number of equality and/or inequality constraints.
The solver has been demonstrated to be fairly robust and highly efficient in solving
various classes of SDP problems arising from the relaxation of
combinatorial optimization problems such as
maximum stable set problems, quadratic assignment problems,
frequency assignment problems, binary quadratic integer programming
problems. It has also worked well on SDP problems arising
from the relaxation of robust clustering problems, rank-one tensor
approximation problems, as well as problems arising
from electronic structure calculations in quantum chemistry.

Our solver is expected to work well on nondegenerate well-posed SDP problems,
but much more future work must be done to make the solver to work well on
degenerate and/or ill-posed problems.
Currently our solver is not catered to problems with SOCP or exponential cone constraints.
As an obvious extension, we are currently extending the
solver to handle problems with the aforementioned cone constraints.

We have also designed a basic user friendly interface for the
user to input their SDP model into the solver.
One of our future works is to expand the flexibility and capability of the
interface such as the ability to handle Hermitian matrices.


\end{document}